\documentclass[12pt,a4paper]{article}

\pdfoutput=1
\usepackage{graphicx}
\usepackage{latexsym}
\usepackage{amssymb}
\usepackage{amsmath,color}
\usepackage{amsthm}
\usepackage{amscd}
\usepackage[dvipsnames]{xcolor}
\usepackage{epstopdf}
\usepackage{amsmath}
\usepackage[all,cmtip]{xy}
\usepackage [latin1]{inputenc}
\newtheorem{Th}{Theorem}[section]
\newtheorem{Def}[Th]{Definition}

\newtheorem{co}[Th]{Corollary}

\newtheorem{pro}[Th]{Proposition}
\newtheorem{re}[Th]{Remark}

\def\be#1\ee{\begin{equation}#1\end{equation}}
\def\bay#1\eay{\!\!\!\left\{\!\!\begin{array}{l}#1\displaystyle\end{array}\right.}
\def\bln#1\eln{\begin{aligned}#1\end{aligned}}
\def\bma#1\ema{{\allowdisplaybreaks\begin{align}#1\end{align}}}

\def\bgr#1\egr{{\allowdisplaybreaks\begin{gather}#1\end{gather}}}
\def\ef#1{(\ref{#1})}

\newcommand{\enD}{\hfill $\Box$\vspace{3truemm} \par}
\newcommand{\bx}{\mbox{\boldmath $x$}}
\newcommand{\bX}{\mbox{\boldmath $X$}}
\newcommand{\e}{\mbox{\boldmath $e$}}

\newcommand{\bv}{\mbox{\boldmath $v$}}

\newcommand{\by}{\mbox{\boldmath $y$}}

\newcommand{\bw}{\mbox{\boldmath $w$}}

\newcommand{\bn}{\mbox{\boldmath $n$}}

\newcommand{\g}{\mbox{\boldmath $\gamma$}}

\newcommand{\m}{\mbox{\boldmath $\mu$}}

\begin{document}

\begin{sloppypar}
\title{\bf Focal surfaces and evolutes of framed curves in hyperbolic 3-space from the viewpoint of Legendrian duality  }
\author{Haibo Yu, Liang Chen$^{*}$ \\
{\scriptsize{\it  School of Mathematics, Changchun Normal University, Changchun 130024, P.R.CHINA
}}\\
{\scriptsize{\it Email address: yuhaibo@ccsfu.edu.cn }}\\
{\scriptsize{\it School of Mathematics and Statistics, Northeast Normal University, Changchun 130024, P.R.CHINA\/}}\\
{\scriptsize{\it Email address:  chenl234@nenu.edu.cn (\scriptsize $^*$ Corresponding author)}}}

\date{ }
\maketitle
\par
\begin{abstract}
A hyperbolic framed curve is a smooth curve with a moving frame in hyperbolic 3-space. It may have singularities. By using this moving frame, we can investigate the differential geometry properties of curves, even at singular points. In fact, we consider the focal surfaces and evolutes of hyperbolic framed curves by using Legendrian dualities which developed by Chen and Izumiya. The focal surfaces are the dual surfaces of tangent indicatrix of original curves. Moreover, classifications of singularities of the serval dual surfaces are shown. By this, we give the relationship among focal surfaces, evolutes and dual surfaces of evolutes. Finally, we study duality relations of singularities between focal surfaces and dual surfaces of evolutes.

\end{abstract}
\renewcommand{\thefootnote}{\fnsymbol{footnote}}

\section{Introduction}
\par
 A mandala of Legendrian dualities for pseudo-spheres in Lorentz-Minkowski space was shown by Izumiya in \cite{15}. It is a generalization of the classical projective duality and the spherical duality. Moreover, Izumiya and the second author of this paper have generalized it into pseudo-spheres in semi-Euclidean space with general index in \cite{3}. It is now a fundamental tool for study of extrinsic differential geometries on submanifolds in these pseudo-spheres from view point of singularity theory (cf. \cite{14,15,16,17} ). We will use it to investigate submanifolds in hyperbolic 3-space and de Sitter 3-space from the viewpoint of Legendrian duality.

The evolutes and focal surfaces of regular curves in different spaces are classical and vital
objects (cf. \cite{2,4,6,8,10,11,16,17}). As we know, evolute of plane curve in $\Bbb R^2$ can be consider as the envelope of the family of normal lines of it in \cite{2,6,7}. Moreover, as the increase of dimension, the envelope of the family of the normal planes of the space curve in $\Bbb R^3$ is surface which is called the focal surface of the original curve in \cite{10,11}. Similarly, the spherical evolute of curves in unit sphere is consider as the envelope of the family of normal geodesics of the original curve. Meanwhile, in \cite{12}, it is also the dual curve of the unit tangent vector of the original curve from the the spherical duality viewpoint. Moreover, for a curve in the unit 3-sphere, the dual of the tangent indicatrix of a curve is naturally surface. It is called the focal surface of the original curves in \cite{17}.  This result gives us the motivation to define focal surfaces of regular curves in pseudo spheres by using Legendrian duality and to investigate some duality properties. Furthermore, in \cite{12}, the evolutes of the original curve in unit 3-sphere can be consider as the critical locus of the focal surface is the evolute of the original curve.

On the other hand, a framed curve in Euclidean 3-space is a smooth curve with a moving frame in \cite{9}. It is an important class of singular curves which may admit singularities.  Moreover, Evolutes and focal surfaces of framed curves in Euclidean 3-space are investigated by Honda and Takahashi in \cite{10}. Furthermore, Hayashi, Izumiya and et al investigate evolutes and focal surfaces of regular curves in hyperbolic 3-space (cf. \cite{13}). Unfortunately, there are not any result in the study of the evolutes and focal surfaces of the curves with singularities in hyperbolic 3-space. In fact, the evolute and focal surface may have singularities, even if the space curve is regular. Therefore, it is essential to investigate the evolutes and focal surfaces of curves with singularities in hyperbolic 3-space to reflect local geometry properties not only at regular points but also at singular points.

After introducing the concept of hyperbolic framed curves in \cite{22}, in this paper, we consider the focal surfaces and evolutes of hyperbolic framed curves by using the Legendrian duality. In Section 2, we introduce the the concepts about $\Delta_{i}$-dual $(i=1,5)$ and frontal. Besides, we give a brief review of definition and properties of framed curves in hyperbolic 3-space. In Section 3, we define the focal surfaces of hyperbolic framed curves. The focal surfaces have two types. The hyperbolic focal surfaces and de sitter focal surfaces can be consider as the $\Delta_{1}$-dual surfaces and $\Delta_{5}$-dual surfaces of unit tangent vector of hyperbolic framed curves respectively. In Section 4, by the criterions of singularities of frontal or front, we give the classifications of singularities of focal surfaces. This reveals the relationship between singular points and invariants of the hyperbolic framed curves. In Section 5, we define two types of evolutes of hyperbolic framed curves and study the relationship between focal surfaces and evolutes from the viewpoint of singularity theory. Moreover, by using Legendrian dualities, we also investigate the singularities of $\Delta_{i}$-dual surfaces of evolutes. In Section 6, we consider the dual relationship of singularities between focal surfaces
and $\Delta_{i}$-dual surfaces of evolutes.

\section{Preliminary}
In this section, we give some basic notions about Legendrian duality, frontal in $H^3$ or $S^3_1$and hyperbolic framed curves. Firstly, we give the definitions of Legendrian duality and frontal which are used in this paper.

\subsection{ Legendrian duality and frontal}

Let $\Bbb R^4=\{(x_0,x_1,x_2,x_3)|x_{i}\in\Bbb R (i=0,1,2,3)\}$ be an 4-dimensional vector space. For any $\bx=(x_0,x_1,x_2,x_3)$, $\by=(y_0,y_1,y_2,y_3)\in\Bbb R^4$, the {\it pseudo scalar product} of $\bx$ and $\by$ can be defined by
$$\langle \bx,\by\rangle=-x_0y_0+\sum_{i=1}^{i=3}x_{i}y_{i},$$
We call $(\Bbb R^4,\langle,\rangle)$ {\it Minkowski space}. We denote it as $\Bbb R^4_1$. A non-zero vector $\bx\in \Bbb R^4_1$ is {\it spacelike, lightlike, or timelike} if $\langle \bx,\bx\rangle>0$, $\langle \bx,\bx\rangle=0$ or $\langle \bx,\bx\rangle<0$, respectively.

  For any $\bx_1,\bx_2,\bx_3\in\Bbb R^4_1$, we define $\bx_1\wedge\bx_2\wedge\bx_3$ by
  $$\bx_1\wedge\bx_2\wedge\bx_3={\rm det}\begin{pmatrix}
   -\e_0&\e_1&\e_2&\e_3\\
   x_0^1&x_1^1&x_2^1&x_3^1\\
   x_0^2&x_1^2&x_2^2&x_3^2\\
   x_0^3&x_1^3&x_2^3&x_3^3\\
   \end{pmatrix}$$
 where $\e_0, \e_1, \e_2, \e_3$ is the canonical basis of $\Bbb R^4_1$. We can easily see that $\langle \bx_{0},\bx_1\wedge\bx_2\wedge\bx_3\rangle={\rm det}(\bx_{0},\bx_1,\bx_2,\bx_3)$ so that $\bx_1\wedge\bx_2\wedge\bx_3$ is pseudo orthogonal to any $x_{i}$, $i=1,2,3$.

We define hyperbolic 3-space by
 $$ H^3=\{\bx\in\Bbb R^4_1|\langle \bx,\by\rangle=-1\},$$
 and de Sitter 3-space
  $$ S^3_1=\{\bx\in\Bbb R^4_1|\langle \bx,\by\rangle=1\}.$$
Then we now consider the following two double fibrations (cf. \cite{3,15})
\begin{align*}
{\rm (1)}~&{\rm (a)}~ H^{n}_{r-1}\times S^{n}_{r}\supset\Delta_{1}=\{(\bv,\bw)|\langle \bv,\bw\rangle=0\},\\
 &{\rm (b)}~ \pi_{11}:\Delta_{1}\to H^{n}_{r-1}, ~~ \pi_{12}:\Delta_{1}\to S^{n}_{r},\\
&{\rm (c)}~ \theta_{11}=\langle d\bv,\bw\rangle|_{\Delta_{1}},~~\theta_{12}=\langle\bv, d\bw\rangle|_{\Delta_{1}},\\
{\rm (2)}~&{\rm (a)}~ S^{n}_{r}\times S^{n}_{r}\supset\Delta_{5}=\{(\bv,\bw)|\langle \bv,\bw\rangle=0\},\\
 &{\rm (b)}~ \pi_{51}:\Delta_{5}\to S^{n}_{r}, ~~ \pi_{52}:\Delta_{5}\to S^{n}_{r},\\
&{\rm (c)}~ \theta_{51}=\langle d\bv,\bw\rangle|_{\Delta_{5}},~~\theta_{52}=\langle\bv, d\bw\rangle|_{\Delta_{5}}.\\
\end{align*}
Here, $\pi_{i1}(\bv,\bw)=\bv$, $\pi_{i2}(\bv,\bw)=\bw$, $\langle d\bv,\bw\rangle=-w_{0}dv_{0}+\sum_{i=1}^{i=3}w_{i}dv_{i}$,
$\langle \bv,d\bw\rangle=-v_{0}dw_{0}+\sum_{i=1}^{i=3}v_{i}dw_{i}$ for $\bv=(v_0,v_1,v_2,v_3)$, $\bw=(w_0,w_1,w_2,w_3)\in \Bbb R^{4}_{1}$.

We remark that $\theta_{i1}^{-1}(0)$ and $\theta_{i2}^{-1}(0)$ define the same tangent hyperplane field over $\Delta_{i}$ which is denoted by $K_{i}$. $(\Delta_{i},K_{i}) (i=1,5)$ is a contact manifold with the contact form $\theta_{ij}$ and both of $\pi_{ij}$ are Legendrian fibrations (see \cite{3,15}).

Let $U\subset\Bbb R^2$ be an open subset. A map $\mathcal{L}_{i}:U\rightarrow\Delta_{i}$ is isotropic if the pull-back $\mathcal{L}^{*}\theta_{ij} (i=1,5,j=1,2)$ vanishes identically (cf. \cite{14}). We say that $\pi_{i1}\circ \mathcal{L}_{i}$ and  $\pi_{i2}\circ \mathcal{L}_{i}$ are $\Delta_{i}$-dual each other if $\mathcal{L}_{i}:U\rightarrow\Delta_{i}$ is isotropic mapping with respect to $K_{i}$ (cf. \cite{14}). This relation is called the Legendrian duality in \cite{15}. Then $f:U\rightarrow H^3$ is said to be a frontal if there a $\Delta_{1}$-dual $g:U\rightarrow S^3_1$. Moreover, if $(f,g):U\rightarrow \Delta_{1}$ is an immersion, we say that $f:U\rightarrow H^3$ is a front (cf. \cite{18}). Similarly, by $\Delta_{5}$-dual, we also can define frontal or front in $S^3_1$.

\subsection{Framed curves in hyperbolic 3-space}
In this subsection, we give a brief review of framed curves in hyperbolic 3-space (cf. \cite{22}).

\begin{Def}
We set that $\g:I\rightarrow H^3$ is a smooth curve and $\Delta_{5}=\{\bv=(\bv_1,\bv_2)\in S^3_1\times S^3_1|\langle\bv_1,\bv_2\rangle=0\}$, then
$(\g,\bv_1,\bv_2):I\rightarrow H^3\times\Delta_{5}$ is a {\it framed curve} if $\langle\g,\bv_{i}\rangle=0$ and $\langle \g',\bv_{i}\rangle=0$ (i=1, 2. $'= \frac{d}{dt}$). We also say that $\g :I\rightarrow H^3$ is a {\it framed curve} (or a {\it framed base curve}) if there exists $\bv:I\rightarrow \Delta_{5}$ such that $(\g,\bv_1,\bv_2):I\rightarrow H^3\times\Delta_{5}$ is a framed curve.
\end{Def}

Since $\bv\in\Delta_{5}$, we put $\m=\g\wedge\bv_1\wedge\bv_2 \in S^3_1$. Then we obtain a moving frame $\{\g,\bv_1,\bv_2,\m\}$ along $\g$ and the following Frenet-Serret type formula
\begin{equation}\label{eq2.1}
\begin{pmatrix}
\g'(t)\\
\bv_1'(t)\\
\bv_2'(t)\\
\m'(t)\\
\end{pmatrix}=
\begin{pmatrix}
0&0&0&m(t)\\
0&0&n(t)&a(t)\\
0&-n(t)&0&b(t)\\
m(t)&-a(t)&-b(t)&0\\
\end{pmatrix}=
\begin{pmatrix}
\g(t)\\
\bv_1(t)\\
\bv_2(t)\\
\m(t)\\
\end{pmatrix},
\end{equation}
where $m(t)=\langle\g'(t),\m(t)\rangle$, $n(t)=\langle\bv_1'(t),\bv_2(t)\rangle$, $a(t)=\langle\bv_1'(t),\m(t)\rangle$, $b(t)=\langle\bv_2'(t),\m(t)\rangle$. And we call the functions $(m,n,a,b)$ the curvature of hyperbolic framed curves.

\begin{Def}
	We say that $(\g,\bv_1,\bv_2)$ and $(\bar{\g},\bar{\bv}_{1},\bar{\bv}_2):I\rightarrow H^3\times \Delta_{5}$ are congruent as hyperbolic framed curves if there exists a matrix $ A\in SO(1,3)$ and a constant vector $\bx \in \Bbb R^4_1$ such that
	$$  \bar{\g}(t)=A(\g(t))+\bx, \ \bar{\bv}_1(t)=A(\bv_1(t)),\ \bar{\bv}_2(t)=A(\bv_2(t)) $$
	for all $t\in I$.
\end{Def}
 Now we give the following the existence and uniqueness for hyperbolic framed curves under the Lorentz motion.
\begin{Th}
	{\bf (Existence theorem for hyperbolic framed curves \cite{22}) } Let $(m,n,a,b):I\rightarrow \Bbb R^4$ be a smooth maping. There exists a hyperbolic framed curve $(\g,\bv_1,\bv_2):I\rightarrow H^3\times \Delta_{5}$ whose the curvature is $(m,n,a,b)$.
\end{Th}

\begin{Th}
{\bf	(Uniqueness of hyperbolic framed curves \cite{22}) } Let $(\g,\bv_1,\bv_2)$ and $(\bar{\g},\bar{\bv}_{1},\bar{\bv}_2):I\rightarrow H^3\times \Delta_{5}$ be hyperbolic framed curves with the curvature $(m,n,a,b)$ and $(\bar{m},\bar{n},\bar{a},\bar{b})$, repectively. Then $(\g,\bv_1,\bv_2)$ and $(\bar{\g},\bar{\bv}_{1},\bar{\bv}_2)$ are congruent as hyperbolic framed curves if and only if $(m,n,a,b)=(\bar{m},\bar{n},\bar{a},\bar{b})$.

\end{Th}

Furthermore, to simplify the expressions in the latter sections, we define four smooth functions as follows:
\begin{equation}
\begin{aligned}
f:I\rightarrow \Bbb R&,\quad f(t)=a(t)b'(t)-a'(t)b(t)+n(t)(a^2(t)+b^2(t)).\\
g:I\rightarrow \Bbb R&,\quad g(t)=m(t)b'(t)-m'(t)b(t)+m(t)a(t)n(t). \\
h:I\rightarrow \Bbb R&,\quad h(t)=m(t)a'(t)-m'(t)a(t)-m(t)b(t)n(t).\\
\sigma:I\rightarrow \Bbb R&,\quad \sigma(t)=f^2(t)-g^2(t)-h^2(t).\\
\end{aligned}
\end{equation}
Now we consider a special moving frame along $\g$ in $H^3$. If $a^2(t)+b^2(t)\neq0$ for all $t \in I$, we put
$$\bn_1(t)=\frac{1}{\sqrt{a^2(t)+b^2(t)}}(a(t)\bv_1(t)+b(t)\bv_2(t)), \quad \bn_2(t)=\frac{1}{\sqrt{a^2(t)+b^2(t)}}(-b(t)\bv_1(t)+a(t)\bv_2(t)).$$
By a direct calculation, $\g\wedge\bv_1\wedge\bv_2=\m$ which is called {\it Frenet type frame} along $\g$ holds and we can verify that $(\g,\bn_1,\bn_2):I\rightarrow H^3\times\Delta_{5}$ is also a framed curve. Then the Frenet-Serret type formula is as follows:
\begin{equation}\label{eq2.3}
\begin{pmatrix}
\g'(t)\\
\bn_1'(t)\\
\bn_2'(t)\\
\m(t)\\
\end{pmatrix}=
\begin{pmatrix}
0&0&0&M(t)\\
0&0&N(t)&A(t)\\
0&-N(t)&0&B(t)\\
M(t)&-A(t)&-B(t)&0\\
\end{pmatrix}
\begin{pmatrix}
\g(t)\\
\bn_1(t)\\
\bn_2(t)\\
\m(t)\\
\end{pmatrix},
\end{equation}
where $M(t)=m(t), N(t)=\frac{f(t)}{a^2(t)+b^2(t)}, A(t)=\sqrt{a^2(t)+b^2(t)}, B(t)=0.$
 When $\bv_1=\bn_1$, $\bv_2=\bn_2$, we calculate that $\sigma=A^2N^2(A^2-M^2)-(MA'-M'A)^2$ under the $\{\g,\bn_1,\bn_2,\m\}$. In this case, we denote $\sigma$ by $\sigma_{F}$.

\section{ Focal surfaces of hyperbolic framed curves}
In this section, we consider the focal surfaces of framed curves $(\g,\bv_1,\bv_2)$ by using the Legendrian dualities for pseudo-spheres in Lorentz-Minkowski space.
\subsection{Hyperbolic focal surfaces}
To consider the hyperbolic focal surfaces, we can define $H^{h}:I\times H^3 \rightarrow \Bbb R$ by $H^h(t,\bx)=\langle \bx, \m(t)\rangle$. Then the discriminant set of
$H^{h}$ can be denoted by $$\mathcal{D}_{H^{h}}=\{\bx\in H^3|H^h(t,\bx)=\frac{\partial H^h}{\partial t}(t,\bx)=0 \}.$$
Since $\{\g,\bv_1, \bv_2, \m\}$ is a basis of $\Bbb R^4_1$ along $\g$, we can set $\bx=x_1\g+x_2\bv_1+x_3\bv_2+x_4\m$, where $x_1, x_2, x_3, x_4 \in \Bbb R$. By a straightforward calculation, we obtain that $H^h(t,\bx)=\frac{\partial H^h}{\partial t}(t,\bx)=0$ if and only if $x_4=0, mx_1+ax_2+bx_3=0$. Since $\bx \in H^3$, it follows that the
parametrization of $\mathcal{D}_{H^{h}}$ is as follows: $$ (t,x_1,x_2,x_3)\mapsto x_1\g(t)+x_2\bv_1(t)+x_3\bv_2(t),$$ where $x_1, x_2, x_3 \in \Bbb R$ satisfy that $mx_1+ax_2+bx_3=0$ and $x_1^2-x_2^2-x_3^2=1$. Then we give the following definition.
\begin{Def}\label{de3.1}
Let $(\g,\bv_1, \bv_2):I\to H^3\times\Delta_{5}$ be a framed curve with the framed curvature $(m, n, a, b)$. Define $\mathcal{F}^{h}: I\times J^{h} \rightarrow H^3$ by $$\mathcal{F}^{h}(t,u_1,u_2,u_3)=u_1\g(t)+u_2\bv_1(t)+u_3\bv_2(t),$$ where $J^{h}=\{(u_1,u_2,u_3)\in \Bbb R^3|mu_1+au_2+bu_3=0, u_1^2-u_2^2-u_3^2=1\}$. We call it the {\it hyperbolic focal surfaces} of framed curve $(\g,\bv_1,\bv_2)$.
\end{Def}

From the viewpoint of envelope, the discriminant set is the envelope of the family $\{\bx\in H^3|H^{h}(t,\bx)=0\}_{t\in I}$. Fixed a $t_0\in I$, since the $\{\bx\in H^3|H^{h}(t,\bx)=0\}_{t_0\in I}$ is the intersection of hyperbolic space $H^3$ and normal hyperplane of framed curve $(\g,\bv_1,\bv_2)$, thus $\mathcal{F}^{h}(I\times J^{h})$ is the (hyperbolic) focal set of framed curve $(\g,\bv_1,\bv_2)$.

At least locally, the condition $(m,a,b)\neq(0,0,0)$ holds. When $m\neq0$, $u_1=-\frac{au_2+bu_3}{m}$ holds. Since $u_1^2-u_2^2-u_3^2=1$, it follows that $(a,b)\neq(0,0)$ holds. Thus, we can rewrite hyperbolic focal surface $\mathcal{F}^{h}$ by using the Frenet type frame $\{\g,\bn_1,\bn_2,\m\}$. By definition \ref{de3.1}, under $A^2-M^2>0$, we show that
$\mathcal{F}^{h}:I\times\Bbb R \rightarrow H^3$,
\begin{equation}
\mathcal{F}^{h}(t,\theta)=\frac{\cosh\theta}{\sqrt{A^2(t)-M^2(t)}}(A(t)\g(t)-M(t)\bn_1(t))+\sinh\theta\bn_2(t).
\end{equation}
Moreover, by equation (\ref{eq2.3}), we have
\begin{equation}\label{eq3.3}
\begin{aligned}
(\mathcal{F}^{h})_t&=\cosh\theta\frac{(MA'-AM')(-M)}{(A^2-M^2)^{3/2}}\g+(\cosh\theta\frac{(MA'-AM')(A)}{(A^2-M^2)^{3/2}}-\sinh\theta N)\bn_1\\
&-\cosh\theta \frac{MN}{\sqrt{A^2-M^2}}\bn_2,\\
(\mathcal{F}^{h})_{\theta}&=\frac{\sinh\theta}{\sqrt{A^2(t)-M^2(t)}}(A(t)\g(t)-M(t)\bn_1(t))+\cosh\theta\bn_2(t).\\
\end{aligned}
\end{equation}
Since $\mathcal{F}^{h}$ satisfy (\ref{eq3.3}), we have $(\mathcal{F}^{h},\m)\in \Delta_1$ is isotropic. Therefore, $\mathcal{F}^{h}$ and $\m$ are $\Delta_1$-dual each other.
Furthermore, by the equations (\ref{eq2.3}) and (\ref{eq3.3}), we set
\begin{equation}\label{eq3.4}
\lambda^{h}(t,\theta)={\rm det}(\mathcal{F}^{h},(\mathcal{F}^{h})_t,(\mathcal{F}^{h})_{\theta},\m)(t,\theta)=\frac{\cosh\theta(MA'-AM')(t)-\sinh\theta (AN\sqrt{A^2-M^2})(t)}{A^2(t)-M^2(t)},
\end{equation}
then we calculate that the singular set $S(\mathcal{F}^{h})=\{(t,\theta)|\lambda^{h}(t,\theta)=0\}$.
\subsection{De Sitter focal surfaces}
We now consider the de Sitter focal surfaces. Define $H^{d}:I\times S^3_1\rightarrow \Bbb R$ by $H^{d}(t,\bx)=\langle\bx,\m(t)\rangle$. Then the discriminant set of $H^{d}$ can be given by
$$\mathcal{D}_{H^{d}}=\{\bx\in S^3_1|H^{d}(t,\bx)=\frac{\partial H^{d}}{\partial t}(t,\bx)=0\}$$
Since $\{\g,\bv_1,\bv_2,\m\}$ is a basis of $\Bbb R^4_1$ along $\g$, we can assume $\bx=x_1\g+x_2\bv_1+x_3\bv_2+x_4\m$, where $x_1, x_2, x_3, x_4\in \Bbb R$. By a direct calculation, we obtain that $H^{d}(t,\bx)=\frac{\partial H^{d}}{\partial t}(t,\bx)=0$ if and only if $x_4=0, mx_1+ax_2+bx_3=0$. Since $\bx\in S^3_1$, the parametrization of $\mathcal{D}_{H^{d}}$ is given by
$$ (t,x_1,x_2,x_3) \mapsto x_1\g(t)+x_2\bv_1(t)+x_3\bv_2(t),$$ where $x_1, x_2, x_3$ satisfy that $mx_1+ax_2+bx_3=0$ and $-x_1^2+x_2^2+x_3^2=1$. Then we give the following definition.

\begin{Def}\label{def3.2}
Let $(\g,\bv_1, \bv_2):I\to H^3\times\Delta_{5}$ be a framed curve with the framed curvature $(m, n, a, b)$. Define $\mathcal{F}^d:I\times J^d\rightarrow S^3_1$ by
$$\mathcal{F}^d(t,u_1,u_2,u_3)=u_1\g(t)+u_2\bv_1(t)+u_3\bv_2(t),$$
where $J^d=\{(u_1,u_2,u_3)\in\Bbb R^3|mu_1+au_2+bu_3=0, -u_1^2+u_2^2+u_3^2=1\}$. We call it the {\it de Sitter focal surface} of framed curve $(\g,\bv_1, \bv_2)$.
\end{Def}

From the view of envelope point, the discriminant set is the envelope of the family $\{\bx\in S^3_1|H^{d}(t,\bx)=0\}_{t\in I}$. For a fixed $t_0 \in I$, since the  $\{\bx\in S^3_1|H^{d}(t_0,\bx)=0\}_{t_0\in I}$ is the intersection of de sitter space $S^3_1$ and normal hyperplane of framed curve $(\g,\bv_1,\bv_2)$, thus $\mathcal{F}^d(I\times J^d)$ is the (de Sitter) foacl set of framed curve $(\g,\bv_1,\bv_2)$.

Suppose that $(a, b)\neq(0,0)$. Then we can rewrite the de Sitter focal surface by using the Frenet type frame $\{\g,\bn_1,\bn_2,\m\}$ as follows. By definition \ref{def3.2}, under the condition $M^2-A^2>0$, we can show
$\mathcal{F}^{d}:I\times\Bbb [0, 2\pi] \rightarrow S_1^3$,
\begin{equation}
\mathcal{F}^d(t,\theta)=\cos\theta\frac{1}{\sqrt{M^2(t)-A^2(t)}}(A(t)\g(t)-M(t)\bn_1(t))+\sin\theta\bn_2(t).
\end{equation}
Furthermore, by equation (\ref{eq2.3}), we have
\begin{equation}\label{eq3.8}
\begin{aligned}
(\mathcal{F}^{d})_t&=\cos\theta\frac{(MA'-AM')(M)}{(M^2-A^2)^{3/2}}\g+(\cos\theta\frac{(MA'-AM')(-A)}{(M^2-A^2)^{3/2}}-\sin\theta N)\bn_1\\
&-\cos\theta \frac{MN}{\sqrt{M^2-A^2}}\bn_2,\\
(\mathcal{F}^{d})_{\theta}&=\frac{-\sin\theta}{\sqrt{M^2(t)-A^2(t)}}(A(t)\g(t)-M(t)\bn_1(t))+\cos\theta\bn_2(t).\\
\end{aligned}
\end{equation}
Since $\mathcal{F}^d$ satisfies (\ref{eq3.8}), we have $(\mathcal{F}^d,\m)\in \Delta_5$ is isotropic. Therefore, $\mathcal{F}^d$ and $\m$ are $\Delta_5$-dual each other. By the equations (\ref{eq2.3}) and (\ref{eq3.8}), we set
\begin{equation}\label{eq3.99}
\lambda^{d}(t,\theta)=\det(\mathcal{F}^d,(\mathcal{F}^d)_t,(\mathcal{F}^d)_{\theta},\m)(t,\theta)=\frac{\cos\theta(MA'-AM') (t)-\sin\theta(AN\sqrt{M^2-A^2})(t)}{M^2(t)-A^2(t)},
\end{equation}
then we calculate that singular set $S(\mathcal{F}^d)=\{(t,\theta)|\lambda^d(t,\theta)=0\}$.

\section{Singularities of focal surfaces of framed curves}
In this section, we give the singularity classifications of focal surfaces of framed curves $(\g,\bv_1,\bv_2)$ by the criterions for singularities at frontal (cf. \cite{5,14}).
\subsection{Singularities of hyperbolic focal surfaces $\mathcal{F}^{h}$}
In this subsection, we give singularity classifications theorem of hyperbolic focal surfaces $\mathcal{F}^{h}$ of framed curves $(\g,\bv_1,\bv_2)$. For the existence conditions of singular points of hyperbolic focal surfaces, we can describe it by using framed curvature $(M,N,A,B)$.
\begin{pro}\label{pro3.2}
Let $p_0=(t_0,\theta_0)$ be a point on hyperbolic focal surfaces $\mathcal{F}^{h}$. When $(MA'-AM',AN\sqrt{A^2-M^2})(t_0)=(0,0)$, all points $(t_0,\theta)$ on $\mathcal{F}^{h}$ are singular points. When $(MA'-AM',AN\sqrt{A^2-M^2})(t_0)\neq(0,0)$, all points $(t_0,\theta)$ on $\mathcal{F}^{h}$ are regular points if and only if $\sigma_{F}(t_0)\leq0$.
\end{pro}
\proof
If $(MA'-AM',AN\sqrt{A^2-M^2})(t_0)=(0,0)$, by (\ref{eq3.4}), we can easily see that all points $(t_0,\theta)$ are singular points on $\mathcal{F}^{h}$. When $(MA'-AM',AN\sqrt{A^2-M^2})(t_0)\neq(0,0)$, we assume that there exists a singular point $(t_0,\theta_0)$. Since $\cosh\theta>0$ and $\cosh\theta^2-\sinh\theta^2=1$, we deduce that $(MA'-AM')^2(t_0)=(1-\frac{1}{\sinh^2\theta_0+1})A^2N^2(A^2-M^2)(t_0)$. When $\theta_0=0$, we obtain $(MA'-AM')(t_0)=0$. Since $\sigma_{F}\leq0$, then  $A^2N^2(A^2-M^2)(t_0)=0$, this is a contradiction. When $\theta_0\neq0$, we can have $\sigma_{F}(t_0)>0$. There is a contradiction. Furthermore, by the proof of necessity, the sufficiency is obvious.
\enD

Since $\mathcal{F}^{h}$ is a $\Delta_{1}$-duality of $\m$, then $\mathcal{F}^{h}$ is a frontal in $H^3$. Assume that $p_0=(t_0,\theta_0)$ is a singular point, then it satisfies that $\lambda^{h}(t_0,\theta_0)=0$. As we known, $p_0$ is an {\it non-degenerate} singular point if $d\lambda^{h}(p_0)\neq0$. By the implicit function theorem, the singular set $S(\mathcal{F}^{h})$ is parameterized by a regular curve $\xi^{h}:I\rightarrow U$ in a neighborhood of an non-degenerate singular point $p_0$. Moreover, there exists a non-zero vector field $\eta^{h}$ near $p_0$ satisfying $\langle\eta^{h}\rangle_{\Bbb R}={\rm ker}d\mathcal{F}^{h}$ along $S(\mathcal{F}^{h})$. We call $\eta^{h}$ {\it null vector field}. Then we can show the following theorem.
\begin{Th}\label{Th3.3}
Under the above notation, we assume point $p_0=(t_0,\theta_0)$ is a singular point of hyperbolic focal surfaces $\mathcal{F}^{h}$. Then we have the followings.\vspace{2mm}\\
\noindent
{\rm(a)} If $(MA'-AM',N)(t_0)\neq(0,0)$, then\\
\mbox{}\quad
{\rm(1)} $p_0$ is cuspidal edge if and only if $\theta'(t_0)-\frac{MN}{\sqrt{A^2-M^2}}(t_0)\neq0$.\\
\mbox{}\quad
{\rm(2)} $p_0$ is swallowtail if and only if $\theta'(t_0)-\frac{MN}{\sqrt{A^2-M^2}}(t_0)=0$, $\frac{d}{dt}(\theta'(t_0)-\frac{MN}{\sqrt{A^2-M^2}})(t_0)\neq0$.\\
\mbox{}\quad
{\rm(3)} $p_0$ is never cuspidal cross cap.\vspace{2mm}\\
\noindent
{\rm(b)} If $(MA'-AM',N)(t_0)=(0,0)$, then\\
\mbox{}\quad
{\rm(1)} $p_0$ is cuspidal edge if and only if  $\cosh\theta_0(MA'-AM')'(t_0)-\sinh\theta_0(AN\sqrt{A^2-M^2})'(t_0)\neq0$\\
\mbox{}\quad
{\rm(2)} $p_0$ is cuspidal beaks if and only if $$\cosh\theta_0(MA'-AM')'(t_0)-\sinh\theta_0(AN\sqrt{A^2-M^2})'(t_0)=0,$$ $$\sinh\theta_0(MA'-AM')'(t_0)-\cosh\theta_0(AN\sqrt{A^2-M^2})'(t_0)\neq0,$$
\begin{align*}
&(\cosh\theta_0(MA'-AM')''(t_0)-\sinh\theta_0(AN\sqrt{A^2-M^2})''(t_0))\sqrt{A^2-M^2}(t_0)\\
&+2M(t_0)N(t_0)(\sinh\theta(MA'-AM')'(t_0)-\cosh\theta(AN\sqrt{A^2-M^2})'(t_0))\neq0.
\end{align*}
\mbox{}\quad
{\rm(3)} $p_0$ is never swallowtail, cuspidal lips and  cuspidal cross cap.
\end{Th}
Here, a singular point $p$ of the map-germ $f : (U,p) \rightarrow M^3$ is called a {\it cuspidal edge} (briefly, {\it CE})  if $f$ at $p$ is $\mathcal{A}$-equivalent to the germ $(u_1,u_2) \mapsto (u_1,u_2^2,u_2^3)$ at 0. A singular point $p$ of the map-germ $f:U \rightarrow S^{3}$ is a {\it swallowtail} (briefly, {\it SW}) if $f$ is $\mathcal{A}$-equivalent to $(u,v)\mapsto (u,4v^{3}+2uv,3v^{4}+uv^{2})$ at 0. A singular point $p$ of $f$ is a {\it cuspidal beaks} (briefly, {\it CBK}) if $f$ is $\mathcal{A}$-equivalent to $(u,v)\mapsto(u,2v^3-u^2v,3v^4-u^2v^2)$ at 0. A singular point $p$ of $f$ is a {\it cuspidal lips} (briefly, {\it CL}) if $f$ is $\mathcal{A}$-equivalent to $(u,v)\mapsto(u,2v^3+u^2v,3v^4+u^2v^2)$ at 0. A singular point $p$ of $f$ is a {\it cuspidal cross cap} (briefly, {\it CCR}) if $f$ is $\mathcal{A}$-equivalent to $(u,v)\mapsto(u,v^2,uv^3)$ at 0. We can draw the pictures of these singularities in the follows. About the criteria for CE, SW, CBK, CL CCR, please see {\cite[Theorem 8.1]{14}}, {\cite[Theorem 8.3]{14}} and {\cite[Theorem 8.5]{14}}.
\begin{center}
 \centering
 \begin{minipage}[c]{0.33\textwidth}
  \centering
  \includegraphics[scale=0.20]{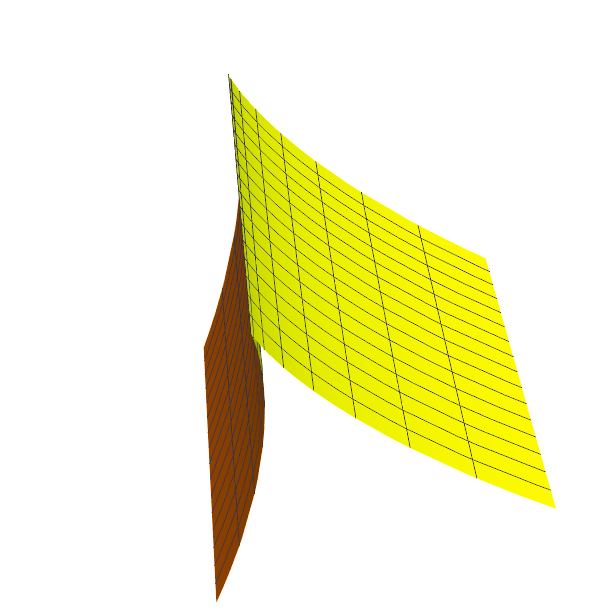}\\
    cuspidal edge
 \end{minipage}%
 \begin{minipage}[c]{0.33\textwidth}
  \centering
  \includegraphics[scale=0.20]{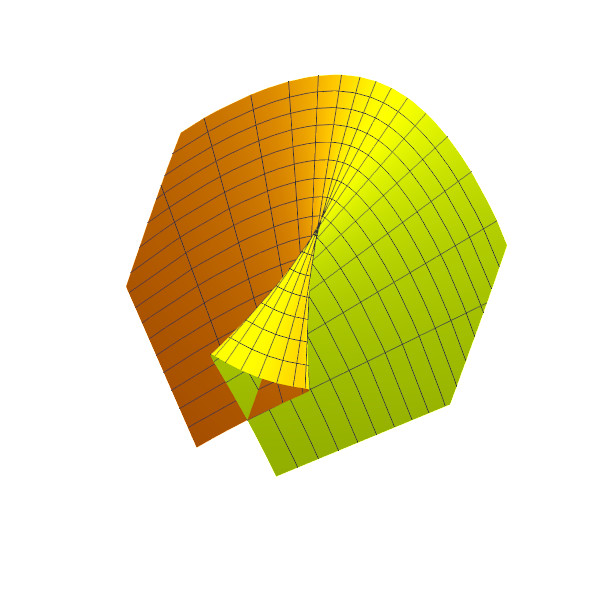}\\
    swallowtail
 \end{minipage}%
 \begin{minipage}[c]{0.33\textwidth}
 \centering
 \includegraphics[scale=0.20]{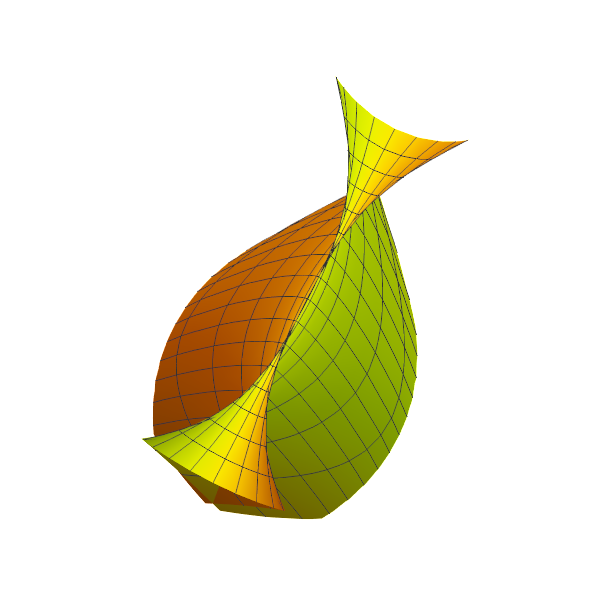}\\
   cuspidal beaks
 \end{minipage}\\
\begin{minipage}[c]{0.3\textwidth}
\centering
\includegraphics[scale=0.3]{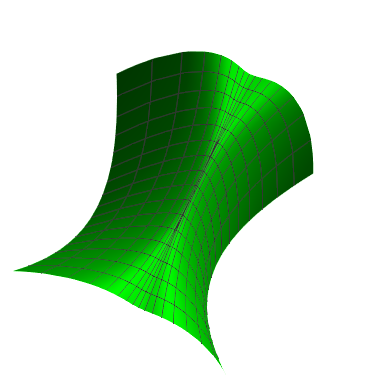}\\
cuspidal lips
\end{minipage}%
\centering
\begin{minipage}[c]{0.3\textwidth}
\centering
\includegraphics[scale=0.33]{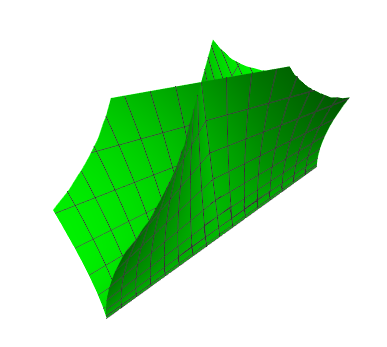}
cuspidal cross cap
\end{minipage}
 \end{center}

\proof
We denote $()_{u}=\partial/\partial u$, $()_{\theta}=\partial /\partial\theta$, $'=\frac{d}{du}$.
By (\ref{eq3.3}), we may take $\eta^{h}=\frac{\partial}{\partial t}+\frac{MN}{\sqrt{A^2-M^2}}\frac{\partial}{\partial \theta}$ on $S(\mathcal{F}^{h})$. Since $(a,b)\neq(0,0)$, $\mathcal{F}^{h}$ at $p_0$ is a front. When $(MA'-AM',N)(t_0)\neq(0,0)$, then $\lambda^{h}(t_0,\theta_0)=0$ and $d\lambda^{h}(t_0,\theta_0)\neq0$ hold. So $p_0$ is non-degenerate. By implicit function theorem, we have singular set $\xi^{h}(t)=(t,\theta(t))$. Then $(\xi^{h})'(t)=\frac{\partial}{\partial t}+\theta'(t)\frac{\partial}{\partial \theta}$ and $\eta^{h}(t)=\frac{\partial}{\partial t}+\frac{MN}{\sqrt{A^2-M^2}}(t)\frac{\partial}{\partial \theta}$ on $\xi^{h}(t)$. We can consider
\begin{equation}
{\rm det}(\eta^{h},(\xi^{h})')(t)={\rm det}\begin{pmatrix}
1&\theta'\\
1&\frac{MN}{\sqrt{A^2-M^2}}\\
\end{pmatrix}(t), \quad
\frac{d}{dt}{\rm det}(\eta^{h},(\xi^{h})')(t)=\frac{d}{dt}(\theta'-\frac{MN}{\sqrt{A^2-M^2}})(t).
\end{equation}
 Thus, ${\rm det}(\eta^{h},(\xi^{h})')(t_0)\neq0$ proves the assertion (1) of (a), ${\rm det}(\eta^{h},(\xi^{h})')(t_0)=0$ and $\frac{d}{dt}{\rm det}(\eta^{h},(\xi^{h})')(t_0)\neq0$ proves the assertion (2) of (a). Since cuspidal cross cap is a frontal which is not a front, then the assertion (3) of (a) holds.

 When $(MA'-AM',N)(t_0)=(0,0)$, we have $\lambda^{h}(t_0,\theta_0)=(\lambda^{h})_{\theta}(t_0,\theta_0)=0$. Thus, $p_0$ is non-degenerate if and only if $\cosh\theta_0(MA'-AM')'(t_0)-\sinh\theta_0(AN\sqrt{A^2-M^2})'(t_0)\neq0$. If $p_0$ is non-degenerate, then singular set $\xi^{h}(\theta)=(t(\theta),\theta)$. Moreover, we have $(\xi^{h})'(\theta)=\partial/\partial \theta$ and $\eta^{h}(t)=\frac{\partial}{\partial t}+\frac{MN}{\sqrt{A^2-M^2}}(t(\theta))\frac{\partial}{\partial \theta}$ on $\xi^{h}(\theta)$. Since
 \begin{equation}
 {\rm det}((\xi^{h})',\eta^{h})(\theta)={\rm det}\begin{pmatrix}
0&1\\
1&\frac{MN}{\sqrt{A^2-M^2}}\\
\end{pmatrix}(\theta)=-1\neq0.
\end{equation}
Therefore, the assertions (1) of (b) holds and there is no swallowtail. In addition, $\mathcal{F}^{h}$ at $p_0$ is also a front in this case. So there is no cuspidal cross cap.

 If $p_0$ is degenerate, we consider
 \begin{equation}
 \det{\rm Hess}\lambda^{h}(p_{0})=\det
\begin{pmatrix}
 (\lambda^{h})_{\theta \theta}&(\lambda^{h})_{\theta t}\\
 (\lambda^{h})_{t\theta}&(\lambda^{h})_{tt}\\
\end{pmatrix}(p_{0}),
 \end{equation}
 Since $(\lambda^{h})_{\theta \theta}(p_0)=\lambda^{h}(p_0)=0$, then $ \det{\rm Hess}\lambda^{h}(p_{0})=-(\lambda^{h})^{2}_{t\theta}(p_0)\leq0$ holds. Thus, there is no cuspidal lips and assertion (3) of (b) holds.
Moreover, we also consider \begin{equation}
 \eta^{h}\eta^{h}\lambda^{h}(p_0)=(\lambda^{h})_{tt}(p_0)+2\frac{MN}{\sqrt{A^2-M^2}}(\lambda^{h})_{t\theta}(p_0).
 \end{equation}
 By a direct calculation, $$(\lambda^{h})_{t\theta}(p_0)=\frac{\sinh\theta_0(A'M-AM')'(t_0)-\cosh\theta(AN\sqrt{A^2-M^2})'(t_0)}{A^2(t_0)-M^2(t_0)}$$ and $$(\lambda^{h})_{tt}(p_0)=\frac{\cosh\theta(A'M-AM')''(t_0)-\sinh\theta(AN\sqrt{A^2-M^2})''(t_0)}{A^2(t_0)-M^2(t_0)},$$ thus $\det{\rm Hess}\lambda^{h}(p_{0})\neq0$ and
 $\eta^{h}\eta^{h}\lambda^{h}(p_0)\neq0$ prove the assertion (2) of (b).
 \enD

\subsection{Singularities of de sitter focal surfaces $\mathcal{F}^{d}$}
In this subsection, we give singularity classifications theorem of of de sitter focal surfaces $\mathcal{F}^d$ of hyperbolic framed curves $(\g,\bv_1,\bv_2)$ under $(a,b)\neq(0,0)$ and $m^2-a^2-b^2>0$. The assertion about the condition for existence of singular points of de sitter $\mathcal{F}^{d}$ is as follows:

\begin{pro}
 When $(MA'-AM',N)(t_0)=(0,0)$, all points $(t_0,\theta)$ on $\mathcal{F}^{d}$ are singular points.
 When $(MA'-AM',N)(t_0)\neq(0,0)$, then there exists a singular point $(t_0,\theta_0)$ on $\mathcal{F}^{d}$.
\end{pro}
\proof
If $(MA'-AM',AN\sqrt{M^2-A^2})(t_0)=(0,0)$, by (\ref{eq3.99}), we can easily see that all points $(t_0,\theta)$ are singular points on $\mathcal{F}^{d}$. When $(MA'-AM',AN\sqrt{M^2-A^2})(t_0)\neq(0,0)$, by (\ref{eq3.99}), the assertion holds.
\enD

Now the singularity classifications theorem of of de sitter focal surfaces $\mathcal{F}^d$ of hyperbolic framed curves $(\g,\bv_1,\bv_2)$ is shown as follows.

\begin{Th}\label{Th4.6}
Under the above notation, we assume point $p_0=(t_0,\theta_0)$ is a singular point of de sitter focal surfaces $\mathcal{F}^{d}$. Then we have the followings.\vspace{2mm}\\
\noindent
{\rm(a)} If $(MA'-AM',N)(t_0)\neq(0,0)$, then\\
\mbox{}\quad
{\rm(1)} $p_0$ is cuspidal edge if and only if $\theta'(t_0)-\frac{MN}{\sqrt{M^2-A^2}}(t_0)\neq0$.\\
\mbox{}\quad
{\rm(2)} $p_0$ is swallowtail if and only if $\theta'(t_0)-\frac{MN}{\sqrt{M^2-A^2}}(t_0)=0$, $\frac{d}{dt}(\theta'(t_0)-\frac{MN}{\sqrt{M^2-A^2}})(t_0)\neq0$.\\
\mbox{}\quad
{\rm(3)} $p_0$ is never cuspidal cross cap.\vspace{2mm}\\
\noindent
{\rm(b)} If $(MA'-AM',N)(t_0)=(0,0)$, then\\
\mbox{}\quad
{\rm(1)} $p_0$ is cuspidal edge if and only if
$\cos\theta_0(MA'-AM')'(t_0)-\sin\theta_0(AN\sqrt{M^2-A^2})'(t_0)\neq0$\\
\mbox{}\quad
{\rm(2)} $p_0$ is cuspidal beaks if and only if
$$\cos\theta_0(MA'-AM')'(t_0)-\sin\theta_0(AN\sqrt{M^2-A^2})'(t_0)=0,$$
$$\sin\theta_0(MA'-AM')'(t_0)+\cos\theta_0(AN\sqrt{M^2-A^2})'(t_0)\neq0,$$
\begin{align*}
&(\cos\theta_0(MA'-AM')''(t_0)-\sin\theta_0(AN\sqrt{M^2-A^2})''(t_0))\sqrt{M^2-A^2}(t_0)\\
&-2M(t_0)N(t_0)(\sin\theta(MA'-AM')'(t_0)+\cos\theta(AN\sqrt{M^2-A^2})'(t_0))\neq0.
\end{align*}
\mbox{}\quad
{\rm(3)} $p_0$ is never swallowtail, cuspidal lips and  cuspidal cross cap.
\end{Th}
\proof  By (\ref{eq3.8}), we have null vector field $\eta^{d}=\frac{\partial}{\partial t}+\frac{MN}{\sqrt{M^2-A^2}}\frac{\partial}{\partial \theta}$ on $S(\mathcal{F}^d)$.
Since $(a,b)\neq0$, then $\mathcal{F}^{d}$ at $p_0$ is a front. When $(MA'-AM',N)(t_0)\neq(0,0)$, $\lambda^{d}(t_0,\theta_0)=0$ and $d\lambda^{d}(t_0,\theta_0)\neq0$. Then $p_0$ is a non-degenerate singular point. By implicit function theorem, we set singular set $\xi^{d}(t)=(t,\theta(t))$.
Then we have $(\xi^{d})'(t)=\frac{\partial}{\partial t}+\theta'(t)\frac{\partial}{\partial \theta}$ and $\eta^{d}(t)=\frac{\partial}{\partial t}+\frac{MN}{\sqrt{M^2-A^2}}(t)\frac{\partial}{\partial \theta}$ on $\xi^{d}(t)$.
We consider
\begin{equation}
{\rm det}(\eta^{d},(\xi^{d})')(t)={\rm det}\begin{pmatrix}
1&\frac{MN}{\sqrt{M^2-A^2}}\\
1&\theta'\\
\end{pmatrix}(t), \quad
\frac{d}{dt}{\rm det}(\eta^{d},(\xi^{d})')(t)=\frac{d}{dt}(\theta'-\frac{MN}{\sqrt{M^2-A^2}})(t),
\end{equation}
 so ${\rm det}(\eta^{d},(\xi^{d})')(t_0)\neq0$ proves the assertion (1) of (a), ${\rm det}(\eta^{d},(\xi^{d})')(t_0)=0$ and $\frac{d}{dt}{\rm det}(\eta^{d},(\xi^{d})')(t_0)\neq0$ proves the assertion (2) of (a). Since cuspidal cross cap is a frontal which is not a front, the assertion (3) of (a) holds.

 When $(MA'-AM',N)(t_0)=(0,0)$, $\lambda^{d}(t_0,\theta_0)=(\lambda^{d})_{\theta}(t_0,\theta_0)=0$. Point $p_0$ is non-degenerate if and only if $\cos\theta_0(MA'-AM')'(t_0)-\sin\theta_0(AN\sqrt{M^2-A^2})'(t_0)\neq0$. Then we have singular set $\xi^{d}(\theta)=(t(\theta),\theta)$. Moreover, we have $(\xi^{d})'(\theta)=\partial/\partial \theta$ and $\eta^{d}(t)=\frac{\partial}{\partial t}+\frac{MN}{\sqrt{M^2-A^2}}(t(\theta))\frac{\partial}{\partial \theta}$ on $\xi^{d}(\theta)$. Since
 \begin{equation}
 {\rm det}((\xi^{d})',\eta^{d})(\theta)={\rm det}\begin{pmatrix}
0&1\\
1&\frac{MN}{\sqrt{A^2-M^2}}\\
\end{pmatrix}(\theta)=-1\neq0.
\end{equation}
Thus, the assertions (1) of (b) holds and there is no swallowtail. In this case, $\mathcal{F}^{d}$ is also a frontal at $p_0$.
  When $p_0$ is degenerate, we consider
 \begin{equation}
 \det{\rm Hess}\lambda^{d}(p_{0})=\det
\begin{pmatrix}
 (\lambda^{d})_{\theta \theta}&(\lambda^{d})_{\theta t}\\
 (\lambda^{d})_{t\theta}&(\lambda^{d})_{tt}\\
\end{pmatrix}(p_{0}),
 \end{equation}
 Since $(\lambda^d)_{\theta \theta}(p_0)=\lambda^{d}(p_0)=0$, then we have that
 $ \det{\rm Hess}\lambda^{d}(p_{0})=-(\lambda^{d})^{2}_{t\theta}(p_0)\leq0$. It follows that there is no cuspidal lips. Therefore,
 the assertion (3) of (b) holds.
Moreover, by a direct calculation, \begin{equation}
 \eta^{d}\eta^{d}\lambda^{d}(p_0)=(\lambda^{d})_{tt}(p_0)+2\frac{MN}{\sqrt{M^2-A^2}}(\lambda^{d})_{t\theta}(p_0).
 \end{equation}
 Since $$(\lambda^{d})_{t\theta}(p_0)=\frac{-\sin\theta_0(MA'-AM')'(t_0)-\cos\theta_0(AN\sqrt{M^2-A^2})'(t_0)}{M^2(t_0)-A^2(t_0)}$$ and $$(\lambda^{d})_{tt}(p_0)=\frac{\cos\theta_0(MA'-AM')''(t_0)-\sin\theta_0(AN\sqrt{M^2-A^2})''(t_0)}{M^2(t_0)-A^2(t_0)},$$ then $\det{\rm Hess}\lambda^{d}(p_{0})\neq0$ and
 $\eta^{d}\eta^{d}\lambda^{d}(p_0)\neq0$ prove the assertion (2) of (b).
 \enD

\section{The properties of evolutes of hyperbolic framed curves}
In this section, we define the evolutes of hyperbolic framed curves and give relationship of focal surfaces and evolutes from the viewpoint of singularity. Moreover, by using the Legendrian dualities, we define the $\Delta_{i}$-dual $(i=1,5)$ surfaces of evolutes and investigate the singularities of it.

\subsection{Relationship between hyperbolic evolutes and hyperbolic focal surfaces}
In this subsection, we give the definition of hyperbolic evolutes of hyperbolic framed curves. Moreover, we give the relationship of hyperbolic focal surfaces and hyperbolic evolutes from the viewpoint of singularity.

We define the secondary discriminant $\mathcal{D}_{H^{h}_{\mu}}^2$ of $H^{h}_{\mu}$, where $$\mathcal{D}_{H^{h}_{\mu}}^2=\{\bx\in H^3|H^h_{\mu}(t,\bx)=\frac{\partial H^h_{\mu}}{\partial t}(t,\bx)=\frac{\partial^2 H^h_{\mu}}{\partial t}(t,\bx)=0 \}.$$

Moreover, by the the secondary discriminant $\mathcal{D}_{H^{h}_{\mu}}^2$, we can give the following the definition of hyperbolic evolutes of hyperbolic framed curves.
\begin{Def}\label{de3.3}
 {\rm The {\it hyperbolic evolute} $\mathcal{E}^{h}_{\g}:I\rightarrow H^3$ of framed curves $(\g,\bv_1,\bv_2):I\rightarrow H^3\times\Delta_{5}$ with $\sigma>0$ is given by}

 $$\mathcal{E}^{h}_{\g}(t)=\frac{1}{\sigma(t)}(f(t)\g(t)-g(t)\bv_1(t)+h(t)\bv_2(t))$$

\end{Def}
 Since $\sigma>0$, it follows that $(a,b)\neq(0,0)$. Then under the Frenet type frame $\{\g,\bn_1,\bn_2,\m\}$, we have $\mathcal{E}^{h}_{\g}:I\rightarrow H^3$,
 \begin{equation}\label{eq3.9}
 \mathcal{E}^{h}_{\g}(t)=\frac{1}{\sqrt{\sigma_{F}(t)}}(A^2(t)N(t)\g(t)-M(t)A(t)N(t)\bn_1(t)+(MA'-M'A)(t)\bn_2(t)),
  \end{equation}

Now we show the relationship between the hyperbolic focal surfaces and hyperbolic evolutes.

 \begin{Th}\label{Th3.7}
Let $(\g,\bv_1,\bv_2):I\rightarrow H^3\times\Delta_{5}$ be a framed curve with $\sigma(t)>0$. $p_0=(t_0,\theta_0)$ is a singular point of $\mathcal{F}^{h}$. Then we have the following:\vspace{2mm}\\
\noindent
{\rm (1)} The image of set of non-degenerate singular point of hyperbolic focal surface $\mathcal{F}^{h}$ coincide with the image of $\mathcal{E}^{h}_{\g}$.\\
\noindent
{\rm (2)} The hyperbolic focal surface $\mathcal{F}^{h}$ at $(t_0,\theta_0)$ is $\mathcal{A}$-equivalent to the cuspidal edge if and only if $\mathcal{E}^{h}_{\g}$ at $t_0$ is a regular point.\\
\noindent
{\rm (3)} The hyperbolic focal surface $\mathcal{F}^{h}$ at $(t_0,\theta_0)$ is $\mathcal{A}$-equivalent to the swallowtail of and only if $\mathcal{E}^{h}_{\g}$ at $t_0$ is locally diffeomorphic to {\rm (2,3,4)}-cusp.
\end{Th}
\proof
Under the Frenet type frame $\{\g,\bn_1,\bn_2,\m\}$, we have

\begin{align}
\mathcal{F}^{h}(t,\theta)=&\frac{\cosh\theta}{\sqrt{A^2(t)-M^2(t)}}(A(t)\g(t)-M(t)\bn_1(t))+\sinh\theta\bn_2(t),\\
\mathcal{E}^{h}_{\g}(t)=&\frac{1}{\sqrt{\sigma_{F}(t)}}(A^2N\g-MAN\bn_1+(MA'-M'A)\bn_2)(t).\label{eq3.111}
\end{align}

 We can set $\cosh\theta(t)=\frac{AN\sqrt{A^2-M^2}}{\sqrt{\sigma_{F}}}(t)$ and $\sinh\theta(t)=\frac{MA'-M'A}{\sqrt{\sigma_{F}}}(t)$ under $\sigma_{F}>0$, where $\theta
:I\rightarrow\Bbb R$.
Then from (\ef{eq3.111}), we can calculate that $$\mathcal{E}^{h}_{\g}(t)=\frac{\cosh\theta(t)}{\sqrt{A^2(t)-M^2(t)}}(A(t)\g(t)-M(t)\bn_1(t))+\sinh\theta(t)\bn_2(t),$$ where $\cosh\theta(t)(MA'-M'A)(t)-\sinh\theta(t)(AN\sqrt{A^2-M^2})(t)=0$. Since $\mathcal{E}^{h}_{\g}$ satisfies (\ref{eq3.4}) and $(MA'-M'A, N)\neq0$, then $\mathcal{E}^{h}_{\g}(t)$ is non-degenerate singular value of $\mathcal{F}^{h}(t,\theta)$. On the other hand, we obtain that the non-degenerate singular set $\xi^{h}(t)$ of $\mathcal{F}^{h}(t,\theta)$ under $(MA'-M'A, N)\neq0$ satisfies $\sigma_{F}(t)>0$ by Proposition \ref{pro3.2} and Theorem \ref{Th3.3}. Thus, we have $\mathcal{F}^{h}(\xi^{h}_{\mu}(t))=\mathcal{E}^{h}_{\g}(t)$. The assertion (1) holds.

To simplify the follow equations, we can define three functions $\varepsilon_1^{h},\varepsilon_2^{h},\varepsilon_3^{h}:I\rightarrow \Bbb R$ as follows:
\begin{equation}\label{eq3.10}
\begin{aligned}
\varepsilon_1^{h}(t)&=\frac{\sinh\theta(t) A(t)}{\sqrt{A^2(t)-M^2(t)}}(\theta'-\frac{MN}{\sqrt{A^2-M^2}})(t),\\
\varepsilon_2^{h}(t)&=\frac{\sinh\theta(t) (-M)(t)}{\sqrt{A^2(t)-M^2(t)}}(\theta'-\frac{MN}{\sqrt{A^2-M^2}})(t),\\
\varepsilon_3^{h}(t)&=\cosh\theta(t)(\theta'-\frac{MN}{\sqrt{A^2-M^2}})(t).\\
\end{aligned}
\end{equation}
By (\ref{eq2.3}) and (\ref{eq3.10}), we obtain that

\begin{align}
(\mathcal{E}_{\g}^{h})'(t)=&\varepsilon_1^{h}(t)\g(t)+\varepsilon_2^{h}(t)\bn_1(t)+\varepsilon_3^{h}(t)\bn_2(t)\label{eq3.11}\\
(\mathcal{E}_{\g}^{h})''(t)=&(\varepsilon_1^{h})'(t)\g(t)+(\varepsilon_2^{h})'(t)\bn_1(t)+(\varepsilon_3^{h})'(t)\bn_2(t)\notag\\
 &+\varepsilon_1^{h}(t)\g'(t)+\varepsilon_2^{h}(t)\bn_1^{'}(t)+\varepsilon_3^{h}(t)\bn_2^{'}(t)\label{eq3.12}\\
 (\mathcal{E}_{\g}^{h})'''(t)=&(\varepsilon_1^{h})''(t)\g(t)+(\varepsilon_2^{h})''(t)\bn_1(t)+(\varepsilon_3^{h})''(t)\bn_2(t)\notag\\
 &+2(\varepsilon_1^{h})'(t)\g'(t)+2(\varepsilon_2^{h})'(t)\bn_1^{'}(t)+2(\varepsilon_3^{h})'(t)\bn_2^{'}(t)\notag\\
 &+\varepsilon_1^{h}(t)\g''(t)+\varepsilon_2^{h}(t)\bn_1^{''}(t)+\varepsilon_3^{h}(t)\bn_2^{''}(t)\label{eq3.13}
\end{align}
By (\ref{eq3.10}) and (\ref{eq3.11}), $(\mathcal{E}_{\g}^{h})'(t)\neq0$ if  and only if $\theta'-\frac{MN}{\sqrt{A^2-M^2}}\neq0$. By the assertion (1) of (a) of Theorem \ref{Th3.3}, the assertion (2) holds.

If $p_0$ is a (2,3,4)-cusp of $\mathcal{E}^{h}_{\mu}$, then $(\mathcal{E}_{\g}^{h})'(t_0)=0$, rank$((\mathcal{E}_{\g}^{h})'',(\mathcal{E}_{\g}^{h})''')(t_0)=2$. $(\mathcal{E}_{\g}^{h})'(t_0)=0$ if and only if $\theta'-\frac{MN}{\sqrt{A^2-M^2}}(t_0)=0$. The necessity of assertion (3) is obvious. Next we consider the sufficiency of it. If $\frac{d}{dt}(\theta'-\frac{MN}{\sqrt{A^2-M^2}})(t_0)\neq0$, we assume that rank$((\mathcal{E}_{\g}^{h})'',(\mathcal{E}_{\g}^{h})''')(t_0)=1$, by (\ref{eq3.12}) and (\ref{eq3.13}), we calculate that
\begin{equation}
\frac{(\varepsilon_1^{h})''}{(\varepsilon_1^{h})'}(t_0)=\frac{(\varepsilon_3^{h})''+2N(\varepsilon_2^{h})'}{(\varepsilon_3^{h})'}(t_0),\quad
\frac{(\varepsilon_2^{h})''-2N(\varepsilon_3^{h})'}{(\varepsilon_2^{h})'}(t_0)=\frac{(\varepsilon_3^{h})''+2N(\varepsilon_2^{h})'}{(\varepsilon_3^{h})'}(t_0).
\end{equation}
Since $\sigma_{F}>0$, it follows that $AN(A^2-M^2)\neq0$. Under the $AN(A^2-M^2)\neq0$, by a long but straightforward computation, we have
$$\frac{-\sigma_{F}}{A^2N(A^2-M^2)}(t_0)=0.$$
This is a contradiction. Thus, by the assertion (2) of (a) of Theorem \ref{Th3.3}, the assertion (3) holds.
\enD

\subsection{Relationship between de sitter evolutes and de sitter focal surfaces}
In this subsection, we give the definition of de sitter evolutes of hyperbolic framed curves. Moreover, we give the relationship of de sitter focal surfaces and de sitter evolutes from the viewpoint of singularity.

We can define de sitter evolutes from the secondary discriminant $\mathcal{D}^2_{H^{d}}$ of $H^{d}$, where $$\mathcal{D}^2_{H^{d}}=\{
\bx\in S^3_1|H^{d}(t,\bx)=\frac{\partial H^d}{\partial t}(t,\bx)=0\}.$$

Moreover, by the secondary discriminant $\mathcal{D}^2_{H^{d}}$, we can give the following the definition of de sitter evolutes

\begin{Def}\label{de3.3}
 {\rm The {\it de sitter evolute} $\mathcal{E}^{d}_{\g}:I\rightarrow S^3_1$ of framed curves $(\g,\bv_1,\bv_2):I\rightarrow H^3\times\Delta_{5}$ with $\sigma<0$ is given by}

 $$\mathcal{E}^{d}_{\g}(t)=\frac{1}{\sqrt{-\sigma(t)}}(f(t)\g(t)-g(t)\bv_1(t)+h(t)\bv_2(t))$$

 \end{Def}
We assume $(a,b)(p_0)\neq0$, under the Frenet type frame $\{\g,\bn_1,\bn_2,\m\}$, by the definition~\ref{de3.3}, under the condition $M^2>A^2$, we have $\mathcal{E}^{d}_{\g}:I\rightarrow S^3_1$,
\begin{equation}\label{eq4.26}
 \mathcal{E}^{d}_{\g}(t)=\frac{1}{\sqrt{-\sigma_{F}}}(A^2N\g-MAN\bn_1+(MA'-M'A)\bn_2),
  \end{equation}

 Then we show the following relationship between the de sitter focal surfaces and de sitter evolutes.
\begin{Th}\label{Th4.8}
Let $(\g,\bv_1,\bv_2):I\rightarrow H^3\times\Delta_{5}$ be a framed curve with $\sigma<0$, $(a,b)\neq(0,0)$ and $m^2-a^2-b^2>0$. $p_0=(t_0,\theta_0)$ is a singular point of $\mathcal{F}^{d}$. Then we have the following:\vspace{2mm}\\
\noindent
{\rm (1)} The image of set of non-degenerate singular point of de sitter focal surface $\mathcal{F}^{d}$ coincide with the image of $\mathcal{E}^{d}_{\g}$.\\
\noindent
{\rm (2)} The de sitter focal surface $\mathcal{F}^{d}$ at $(t_0,\theta_0)$ is $\mathcal{A}$-equivalent to the cuspidal edge if and only if $\mathcal{E}^{d}_{\g}$ at $t_0$ is a regular point.\\
\noindent
{\rm (3)} The de sitter focal surface $\mathcal{F}^{d}$ at $(t_0,\theta_0)$ is $\mathcal{A}$-equivalent to the swallowtail of and only if $\mathcal{E}^{d}_{\g}$ at $t_0$ is locally diffeomorphic to {\rm (2,3,4)}-cusp.
\end{Th}
Under the Frenet type frame $\{\g,\bn_1,\bn_2,\m\}$, we have

\begin{align}
\mathcal{F}^d(t,\theta)=&\cos\theta\frac{1}{\sqrt{M^2(t)-A^2(t)}}(A(t)\g(t)-M(t)\bn_1(t))+\sin\theta\bn_2(t).\\
\mathcal{E}^{d}_{\g}(t)=&\frac{1}{\sqrt{-\sigma_{F}(t)}}(A^2N\g-MAN\bn_1+(MA'-M'A)\bn_2)(t).\label{eq4.28}
\end{align}

From (\ref{eq4.28}), we easily set $\cos\theta(t)=\frac{AN\sqrt{M^2-A^2}}{\sqrt{-\sigma_{F}}}(t)$ and $\sin\theta(t)=\frac{MA'-M'A}{\sqrt{-\sigma_{F}}}(t)$ under $\sigma_{F}<0$.
Then we calculate that $\mathcal{E}^{d}_{\g}(t)=\frac{\cos\theta(t)}{\sqrt{M^2(t)-A^2(t)}}(A(t)\g(t)-M(t)\bn_1(t))+\sin\theta(t)\bn_2(t)$, where $\cos\theta(t)(MA'-M'A)(t)-\sin\theta(t)(AN\sqrt{M^2-A^2})(t)=0$. Since $\mathcal{E}^{d}_{\g}$ satisfies (\ref{eq3.99}) and $(MA'-AM',N)\neq(0,0)$, then $\mathcal{E}^{d}_{\g}(t)$ is non-degenerate singular value of $\mathcal{F}^{d}(t,\theta)$. On the other hand, by $M^2>A^2>0$, we obtain that the only non-degenerate singular set $\xi^{d}(t)$ of $\mathcal{F}^{d}(t,\theta)$ satisfies $\sigma_{F}(t)<0$, thus $\mathcal{F}^{d}(\xi^{d}(t))=\mathcal{E}^{d}_{\g}(t)$. The assertion (1) holds.

To make the expressions of following equations simple, we define three functions $\varepsilon_1^d, \varepsilon_2^d, \varepsilon_3^d:I\rightarrow \Bbb R$ as follows
\begin{equation}\label{eq4.29}
\begin{aligned}
\varepsilon_1^{d}(t)&=\frac{\sin\theta(t) (-A)(t)}{\sqrt{M^2(t)-A^2(t)}}(\theta'-\frac{MN}{\sqrt{M^2-A^2}})(t),\\
\varepsilon_2^{d}(t)&=\frac{\sin\theta(t) M(t)}{\sqrt{M^2(t)-A^2(t)}}(\theta'-\frac{MN}{\sqrt{M^2-A^2}})(t),\\
\varepsilon_3^{d}(t)&=\cos\theta(t)(\theta'-\frac{MN}{\sqrt{M^2-A^2}})(t).\\
\end{aligned}
\end{equation}
By (\ref{eq2.3}) and (\ref{eq4.29}), we obtain that

\begin{align}
(\mathcal{E}_{\g}^{d})'(t)=&\varepsilon_1^{d}(t)\g(t)+\varepsilon_2^{d}(t)\bn_1(t)+\varepsilon_3^{d}(t)\bn_2(t)\label{eq4.30}\\
(\mathcal{E}_{\g}^{d})''(t)=&(\varepsilon_1^{d})'(t)\g(t)+(\varepsilon_2^{d})'(t)\bn_1(t)+(\varepsilon_3^{d})'(t)\bn_2(t)\notag\\
 &+\varepsilon_1^{d}(t)\g'(t)+\varepsilon_2^{d}(t)\bn_1^{'}(t)+\varepsilon_3^{d}(t)\bn_2^{'}(t)\label{eq4.31}\\
 (\mathcal{E}_{\g}^{d})'''(t)=&(\varepsilon_1^{d})''(t)\g(t)+(\varepsilon_2^{d})''(t)\bn_1(t)+(\varepsilon_3^{d})''(t)\bn_2(t)\notag\\
 &+2(\varepsilon_1^{d})'(t)\g'(t)+2(\varepsilon_2^{d})'(t)\bn_1^{'}(t)+2(\varepsilon_3^{d})'(t)\bn_2^{'}(t)\notag\\
 &+\varepsilon_1^{d}(t)\g''(t)+\varepsilon_2^{d}(t)\bn_1^{''}(t)+\varepsilon_3^{d}(t)\bn_2^{''}(t)\label{eq4.32}
\end{align}

By (\ref{eq4.29}) and (\ref{eq4.30}), $(\mathcal{E}_{\g}^{d})'(t)\neq0$ if  and only if $\theta'-\frac{MN}{\sqrt{A^2-M^2}}\neq0$. By the assertion (1) of (a) of Theorem \ref{Th4.6}, the assertion (2) holds.

If $p_0$ is a (2,3,4)-cusp of $\mathcal{E}^{d}_{\g}$, then $(\mathcal{E}_{\g}^{d})'(t_0)=0$, rank$((\mathcal{E}_{\g}^{d})'',(\mathcal{E}_{\g}^{d})''')(t_0)=2$. The necessity is obvious. We consider the sufficiency. $(\mathcal{E}_{\g}^{d})'(t_0)=0$ if and only if $\theta'-\frac{MN}{\sqrt{A^2-M^2}}(t_0)=0$. Moreover, if $\frac{d}{dt}(\theta'-\frac{MN}{\sqrt{A^2-M^2}})(t_0)\neq0$, we can assume that rank$((\mathcal{E}_{\g}^{d})'',(\mathcal{E}_{\g}^{d})''')(t_0)=1$, using (\ref{eq4.31}) and (\ref{eq4.32}), a direct calculation gives that
\begin{equation}\label{eq4.33}
\frac{(\varepsilon_1^{d})''}{(\varepsilon_1^{d})'}(t_0)=\frac{(\varepsilon_3^{d})''+2N(\varepsilon_2^{d})'}{(\varepsilon_3^{d})'}(t_0),\quad
\frac{(\varepsilon_2^{d})''-2N(\varepsilon_3^{d})'}{(\varepsilon_2^{d})'}(t_0)=\frac{(\varepsilon_3^{d})''+2N(\varepsilon_2^{d})'}{(\varepsilon_3^{d})'}(t_0).
\end{equation}
Since $\sigma_{F}<0$, it follows that $(MA'-AM',N)\neq(0,0)$. If $N\neq0$, by a long but straightforward computation, we have
$$\frac{-\sigma_{F}}{A^2N(M^2-A^2)}=0,$$
This is a contradiction. If $MA'-AM'\neq0$, by a long but straightforward computation, we have
$$ \frac{-\sigma_{F}}{A^2}=0,$$
This is a contradiction. Thus, rank$((\mathcal{E}_{\g}^{d})'',(\mathcal{E}_{\g}^{d})''')(t_0)=2$. And by the assertion (2) of (a) of Theorem \ref{Th4.6}, the assertion (3) holds.
\enD

\subsection{Singularities of $\Delta_1$-dual surfaces of hyperbolic evolute $\mathcal{E}^{h}_{\g}$}
In this subsection, we study singularities of the $\Delta_1$-dual surfaces of hyperbolic evolute $\mathcal{E}^{h}_{\g}$.
Now we define $H_{\mathcal{E}^{h}}:I\times S^3_1\rightarrow \Bbb R$ by $H_{\mathcal{E}^{h}}(t,\bx)=\langle\bx, \mathcal{E}^{h}_{\g}(t)\rangle$ under the assumption $\sigma>0$, and denote the discriminant set of $H_{\mathcal{E}^{h}}$ by
$$\mathcal{D}_{\mathcal{E}^{h}}=\{\bx\in S^3_1|H_{\mathcal{E}^{h}}(t,\bx)=\frac{\partial H_{\mathcal{E}^{h}}}{\partial t}(t,\bx)=0\}.$$
Since $\sigma\neq0$, it follows that $(a,b)\neq(0,0)$. Then we deduce that $a^2+b^2-m^2\neq0$ by $\sigma_{F}\neq0$. By a straightforward computation, we find that  $$\{\mathcal{E}^{h}_{\g},\m,\frac{(\mathcal{E}^{h}_{\g})^{'}}{\| (\mathcal{E}^{h}_{\g})^{'}\| },\frac{-m\g+a\bv_1+b\bv_2}{\sqrt{-m^2+a^2+b^2}}\}$$ is a basis of $\Bbb R^4_1$ along $\g$. Then we can assume that $\bx=x_1\mathcal{E}^{h}_{\g}+x_2\m+x_3\frac{(\mathcal{E}^{h}_{\g})^{'}}{\| (\mathcal{E}^{h}_{\g})^{'}\| }+x_4\frac{-m\g+a\bv_1+b\bv_2}{\sqrt{-m^2+a^2+b^2}}$, where $x_1,x_2,x_3,x_4 \in \Bbb R$.  By a direct calculation, $H_{\mathcal{E}^{h}}(t,\bx)=\frac{\partial H_{\mathcal{E}^{h}}}{\partial t}(t,\bx)=0$ if and only if $x_1=x_3=0$. Moreover, since $\bx \in S^3_1$, there exists $\theta \in [0,2\pi)$
such that $x_2=\cos\theta, x_4=\sin\theta.$ It follows that the parametrization of $\mathcal{D}_{\mathcal{E}^{h}}$ is given by $$(t,\theta)\mapsto \cos\theta\m(t)+\sin\theta\frac{-m\g+a\bv_1+b\bv_2}{\sqrt{-m^2+a^2+b^2}} $$ and  is denoted by $F^{d}_{\mathcal{E}^{h}}(t,\theta)$, where $F^{d}_{\mathcal{E}^{h}}:I\times[0,2\pi)\rightarrow S^3_1$. We call it {\it $\Delta_1$-dual surfaces} of $\mathcal{E}^{h}_{\g}$ under the assumption $\sigma>0$. Since $(a,b)\neq(0,0)$, then under the Frenet type frame $\{\g,\bn_1,\bn_2,\m\}$, $$F^{d}_{\mathcal{E}^{h}}(t,\theta)=\cos\theta\m(t)+\sin\theta(\frac{-M(t)\g(t)+A(t)\bn_1(t)}{\sqrt{A^2(t)-M^2(t)}}).$$ Furthermore, by equation (\ref{eq2.3}), we have
\begin{equation}\label{eqHE}
\begin{aligned}
(F^{d}_{\mathcal{E}^{h}})_{t}=&(\cos\theta M+\sin\theta\frac{(MA'-AM')A}{(A^2-M^2)^{3/2}})\g-(\cos\theta A+\sin\theta\frac{(MA'-AM')M}{(A^2-M^2)^{3/2}})\bn_1\\
&+(\frac{\sin\theta AN}{\sqrt{A^2-M^2}})\bn_2
+\sin\theta\sqrt{A^2-M^2}\m\\
(F^{d}_{\mathcal{E}^{h}})_{\theta}=&\cos\theta(\frac{-M\g+A\bn_1}{\sqrt{A^2-M^2}})-\sin\theta\m
\end{aligned}
\end{equation}
By the equations (\ref{eq3.9}) and (\ref{eqHE}), we can define $\lambda_{\mathcal{E}^{h}}:I\times[0,2\pi)\rightarrow \Bbb R$ as follows:
\begin{equation}\label{eq4.18}
\lambda_{\mathcal{E}^{h}}(t,\theta)=\rm{ det}(F^{d}_{\mathcal{E}^{h}},(F^{d}_{\mathcal{E}^{h}})_{t},(F^{d}_{\mathcal{E}^{h}})_{\theta},\mathcal{E}^{h}_{\g})(t,\theta)=-\sin\theta\frac{\sqrt{\sigma_{F} (t)}}{A^2-M^2}(t)
\end{equation}

Then we easily obtain that singular set $S(F^{d}_{\mathcal{E}^{h}})=\{(t,\theta)|\sin\theta=0\}$ under $\sigma_{F}>0$. Moreover, by equation (\ref{eqHE}), $(F^{d}_{\mathcal{E}^{h}},\mathcal{E}^{h}_{\g})\in\Delta_1$ must be isotropic. Thus,  $F^{d}_{\mathcal{E}^{h}}$ and $\mathcal{E}^{h}_{\g}$ are $\Delta_1$-dual each other.

\begin{Th}\label{Th4.1}
Under the above notation, we assume point $p_0=(t_0,\theta_0) \in S(F^{d}_{\mathcal{E}^{h}})$ under the assumption $\sigma_{F}>0$. Then we have the followings.\vspace{2mm}\\
\mbox{}\quad
{\rm(a)} $p_0$ is a non-degenerate singular point of the first kind.\\
\mbox{}\quad
{\rm (1)} $p_0$ is a cuspidal edge if and only if $$\frac{AN\sqrt{A^2-M^2}(MA'-M'A)'-(AN\sqrt{A^2-M^2})'(MA'-M'A)}{(AN)^2(A^2-M^2)-(M'A-MA')^2}(t_0)-\frac{MN}{\sqrt{A^2-M^2}}(t_0)\neq0.$$          \\
\mbox{}\quad
{\rm (2)} $p_0$ is cuspidal cross cap if and only if
$$\frac{AN\sqrt{A^2-M^2}(MA'-M'A)'-(AN\sqrt{A^2-M^2})'(MA'-M'A)}{(AN)^2(A^2-M^2)-(M'A-MA')^2}(t_0)-\frac{MN}{\sqrt{A^2-M^2}}(t_0)=0,$$
$$\frac{d}{dt}(\frac{AN\sqrt{A^2-M^2}(MA'-M'A)'-(AN\sqrt{A^2-M^2})'(MA'-M'A)}{(AN)^2(A^2-M^2)-(M'A-MA')^2}-\frac{MN}{\sqrt{A^2-M^2}})(t_0)\neq0.$$
\end{Th}
\proof
Since $(\lambda_{\mathcal{E}^{h}})_{\theta}\neq0$, by implicit theorem, then we can set non-degenerate singular set $\xi_{\mathcal{E}^{h}}(t)=(t,\theta(t))$. The null vector field $\eta_{\mathcal{E}^{h}}=\frac{\partial}{\partial t}+\sqrt{A^2-M^2}(t)\frac{\partial}{\partial \theta}$ on $\xi_{\mathcal{E}^{h}}(t)$. Since $(\lambda_{\mathcal{E}^{h}})_{t}(t,\theta)=0$ on $\xi_{\mathcal{E}^{h}}(t)$, we have $(\xi_{\mathcal{E}^{h}})'(t)=(1,0)$. Then we consider
\begin{equation}
 {\rm det}(\eta_{\mathcal{E}^{h}},(\xi_{\mathcal{E}^{h}})^{'})(t)={\rm det}\begin{pmatrix}
1&\sqrt{A^2-M^2}\\
1&0\\
\end{pmatrix}(t)=-\sqrt{A^2-M^2}(t)\neq0.
\end{equation}
Thus, $p_0$ is a non-degenerate singular point of the first kind. Furthermore, by a long but straightforward computation, we deduce that $\eta_{\mathcal{E}^{h}}\mathcal{E}^{h}_{\g}\neq0$ at $t_0$ if and only if
$$\frac{AN\sqrt{A^2-M^2}(MA'-M'A)'-(AN\sqrt{A^2-M^2})'(MA'-M'A)}{(AN)^2(A^2-M^2)-(M'A-MA')^2}(t_0)-\frac{MN}{\sqrt{A^2-M^2}}(t_0)\neq0.$$ Since $F^{d}_{\mathcal{E}^{h}}$ and $\mathcal{E}^{h}_{\g}$ are $\Delta_1$-dual each other, then the $\mathcal{E}^{h}_{\g}$ can be a unit normal vector of $F^{d}_{\mathcal{E}^{h}}$ in $S^3_1$.  Then $\eta_{\mathcal{E}^{h}}\mathcal{E}^{h}_{\g}\neq0$ at $p_0$ prove that $F^d_{\mathcal{E}^{h}}$ is a front at $p_0$. Thus, the assertion (1) holds.

Let us define a lift $\omega^{h}:U\rightarrow T^{*}_{F^{d}_{\mathcal{E}^{h}}}S^3_1$ by
$$\omega_{p}^{h}(\bv)=\langle \bv, \mathcal{E}_{\g}^{h}(p)\rangle, \bv\in  T_{F^{d}_{\mathcal{E}^{h}}}S^3_1, $$
where $p  $ is in the neighborhood $U$ of $F^{d}_{\mathcal{E}^{h}}$. We set $\pi\circ\omega^{h}:U\rightarrow F^{d}_{\mathcal{E}^{h}}$, then $(\pi\circ\omega^{h})_{*}(r)=dF^{d}_{\mathcal{E}^{h}}(r)$ for any vector $r\in T_{p}U$. Since $\langle dF^{d}_{\mathcal{E}^{h}}(r),\mathcal{E}_{\g}^{h}(p)\rangle$=0, thus we obtain that $(\pi\circ\omega^{h})_{*}(T_{p}U)\subset\ker\omega_{p}^{h}$. This indicates that $\omega^{h}$ is the admissible lift of $F^{d}_{\mathcal{E}^{h}}$.
 By a direct calculation, we obtain that $\frac{(\mathcal{E}_{\g}^{h})^{'}}{\|(\mathcal{E}_{\g}^{h})^{'}\|}$ is transversal to subspace $dF^{d}_{\mathcal{E}^{h}}(T_{p}U)$ at singular point $p$ and satisfies that $\langle \frac{(\mathcal{E}_{\g}^{h})^{'}}{\|(\mathcal{E}_{\g}^{h})^{'}\|}, \mathcal{E}_{\g}^{h}\rangle=0$. Putting \begin{equation*}
\begin{aligned}
\bX(t)=
\begin{cases}
\frac{(\mathcal{E}_{\g}^{h})^{'}(t)}{\|(\mathcal{E}_{\g}^{h})^{'}(t)\|}, \ \ &{\rm if} ~(\mathcal{E}_{\g}^{h})^{'}(t)\neq0,\\
{\lim\limits_{t\rightarrow +t_0}}\frac{(\mathcal{E}_{\g}^{h})^{'}(t)}{\|(\mathcal{E}_{\g}^{h})^{'}\|(t)}, \ \ &{\rm if}~ (\mathcal{E}_{\g}^{h})^{'}(t_0)=0,
\end{cases}
\end{aligned}
\end{equation*}
then $\bX$ ia a non-zero vector field along $F^{d}_{\mathcal{E}^{h}}|_{S(F^{d}_{\mathcal{E}^{h}})}$. Then,
we can set $\psi^{h}=\langle \eta_{\mathcal{E}^{h}}\bX, \mathcal{E}_{\g}^{h}\rangle$ on singular set $S(F^{d}_{\mathcal{E}^{h}})$.
 Since $\langle\bX(t), \mathcal{E}_{\g}^{h}(t)\rangle=0$, we have $\psi^{h}=-\langle \bX, \eta_{\mathcal{E}^{h}}\mathcal{E}^{h}_{\g}\rangle(t)=-\langle \bX,(\mathcal{E}_{\g}^{h})^{'}\rangle(t)$.
To make the following statement clear, we aseume that $(\mathcal{E}_{\g}^{h})^{'}(0)=0$ at $t=0$. By the definition of $\bX$,
\begin{equation}
\begin{aligned}
\psi^{h}(t)=
\begin{cases}
-\langle\frac{(\mathcal{E}_{\g}^{h})^{'}(t)}{\|\mathcal{E}_{\g}^{h})^{'}\|(t)},(\mathcal{E}_{\g}^{h})^{'}(t)\rangle, \ \ &{\rm if} ~(\mathcal{E}_{\g}^{h})^{'}(t)\neq0,\\
-\langle{\lim\limits_{t\rightarrow +0}}\frac{(\mathcal{E}_{\g}^{h})^{'}(t)}{\|\mathcal{E}_{\g}^{h})^{'}\|(t)},(\mathcal{E}_{\g}^{h})^{'}(0)\rangle, \ \ &{\rm if}~ (\mathcal{E}_{\g}^{h})^{'}(0)=0
\end{cases}
\end{aligned}
\end{equation}
Since $(\mathcal{E}_{\g}^{h})^{'}(0)=0$, we can find a deleted neighborhood of $t=0$ such that $(\mathcal{E}_{\g}^{h})^{'}(t)\neq0$. We can deduce that $\psi^{h}(t)=-\frac{(\varepsilon^h)^2(t)}{|\varepsilon^{h}(t)|}\neq0$ in this deleted neighborhood, where

\begin{equation*}
\begin{aligned}
&\varepsilon^h=\frac{AN\sqrt{A^2-M^2}(MA'-AM')'-(AN\sqrt{A^2-M^2})'(MA'-AM')}{A^2N^2(A^2-M^2)-(MA'-M'A)^2}\\
&-\frac{MN}{\sqrt{A^2-M^2}}.\\
\end{aligned}
\end{equation*}
satisfies $\varepsilon^h(t)=0$ if and only if $(\mathcal{E}_{\g}^{h})^{'}(t)=0$.
 Furthermore, we can assume that $\varepsilon^h(t)>0$ on the open interval $(0,t)$,
  by L'Hopital's rule, $\psi^{h}(0)=-{\lim\limits_{t\rightarrow +0}}\frac{(\varepsilon^h)^2(t)}{|\varepsilon^h(t)|}=-2\varepsilon^h(0)=0$.
 Since $\frac{d}{dt}\psi^{h}(t)=-2(\varepsilon^h)'(t)$ on the open interval $(0,t)$, then $\frac{d}{dt}\psi^{h}(0)=-{\lim\limits_{t\rightarrow +0}}2(\varepsilon^h)'(t)=-2(\varepsilon^h)'(0)$.

 By the above, $ \psi^{h}(t)=0$ if and only if $t=0$ if and only if $\varepsilon^h(t)=0$. Thus, for a singular point $p_0=(t_0,\theta(t_0))$, $ \psi^{h}(t_0)=0$ and $ \frac{d}{dt}\psi^{h
}(t_0)\neq0$ if and only if $\varepsilon^h(t_0)=0$ and $(\varepsilon^h)'(t_0)\neq0$. Then $ \psi^{h}(t_0)=0$ and $ \frac{d}{dt}\psi^{h}(t_0)\neq0$ prove that the assertion (2) holds.

\enD
Now we give the relationship between the $\Delta_1$-dual surfaces $F^{d}_{\mathcal{E}^{h}}$ and $\mathcal{E}^{h}_{\g}$ from singular viewpoint.
\begin{pro}\label{pro5.3}
 We assume point $p_0=(t_0,\theta_0) \in S(F^{d}_{\mathcal{E}^{h}})$, then \vspace{2mm}\\
\mbox{}\quad
{\rm (1)} $p_0$ is a cuspidal edge if and only if  $\mathcal{E}^{h}_{\g}$ at $t_0$ is a regular point.\\
\mbox{}\quad
{\rm (2)} $p_0$ is cuspidal cross cap if and only if $\mathcal{E}^{h}_{\g}$ is {\rm (2,3,4)}-cusp at $t_0$.
 \end{pro}
 \proof By the proof of Theorem \ref{Th4.1}, we can easily observe that $p_0$ is a cuspidal edge if and only if $(\mathcal{E}^{h}_{\g})'(t_0)\neq0$ and $p_0$ is a cuspidal cross cap if and only if $\varepsilon^h(t_0)=0$ and $(\varepsilon^h)'(t_0)\neq0$. From the proof of Theorem \ref{Th3.7}, we see that $\varepsilon^h(t_0)=0$ and $(\varepsilon^h)'(t_0)\neq0$ if and if $\mathcal{E}^{h}_{\g}$ is {\rm (2,3,4)}-cusp at $t_0$.
 \enD
By comparing the Theorem \ref{Th3.7} with Proposition \ref{pro5.3}, we have the following corollary about the relationship between singular points of hyperbolic focal surfaces $\mathcal{F}^{h}$ and $\Delta_1$-dual surfaces $F^{d}_{\mathcal{E}^{h}}$.
 \begin{co}
  We assume point $p_0=(t_0,\theta_0)\in S(\mathcal{F}^{h})$, then there exists point $q_0=(t_0,\bar{\theta_0}) \in S(F^{d}_{\mathcal{E}^{h}})$ so that \vspace{2mm}\\
\mbox{}\quad
{\rm (1)} $\mathcal{F}^{h}$ at $(t_0,\theta_0)$ is cuspidal edge if and only if $F^{d}_{\mathcal{E}^{h}}$ at $(t_0,\bar{\theta_0})$ is a cuspidal edge.\\
 \mbox{}\quad
  {\rm (2)} $\mathcal{F}^{h}$ at $(t_0,\theta_0)$ is swallowtail if and only if $F^{d}_{\mathcal{E}^{h}}$ at $(t_0,\bar{\theta_0})$ is a cuspidal cross cap.
 \end{co}
\subsection{Singularities of $\Delta_5$-dual surfaces of de sitter evolute $\mathcal{E}^{d}_{\g}$}
In this subsection, we investigate singularities of the $\Delta_5$-dual surfaces of de sitter evolute $\mathcal{E}^{d}_{\g}$.
Now we define $H_{\mathcal{E}^{d}}:I\times S^3_1\rightarrow \Bbb R$ by $H_{\mathcal{E}^{d}}(t,\bx)=\langle\bx, \mathcal{E}^{d}_{\g}(t)\rangle$ under the assumption $\sigma<0$, and denote the discriminant set of $H_{\mathcal{E}^{d}}$ by
$$\mathcal{D}_{\mathcal{E}^{d}}=\{\bx\in S^3_1|H_{\mathcal{E}^{d}}(t,\bx)=\frac{\partial H_{\mathcal{E}^{d}}}{\partial t}(t,\bx)=0\}.$$

Under $m^2>a^2+b^2$, by a direct computation, we find that  $$\{\mathcal{E}^{d}_{\g},\m,\frac{(\mathcal{E}^{d}_{\g})^{'}}{\| (\mathcal{E}^{d}_{\g})^{'}\| },\frac{-m\g+a\bv_1+b\bv_2}{\sqrt{m^2-a^2-b^2}}\}$$ is a basis along $\g$ in $\Bbb R^4_1$. Set $\bx=x_1\mathcal{E}^{d}_{\g}+x_2\m+x_3\frac{(\mathcal{E}^{d}_{\g})^{'}}{\| (\mathcal{E}^{d}_{\g})^{'}\| }+x_4\frac{-m\g+a\bv_1+b\bv_2}{\sqrt{m^2-a^2-b^2}}$, where $x_1,x_2,x_3,x_4 \in \Bbb R$. Then $H_{\mathcal{E}^{d}}(t,\bx)=\frac{\partial H_{\mathcal{E}^{d}}}{\partial t}(t,\bx)=0$ if and only if $x_1=x_3=0$. Moreover, since $\bx \in S^3_1$, there exists $\theta\in \Bbb R$ such that $x_2=\cosh\theta$, $x_4=\sinh\theta$. It follows that the parametrization of $\mathcal{D}_{\mathcal{E}^{d}}$ is given by $$(t,\theta)\mapsto \cosh\theta\m(t)+\sinh\theta\frac{-m\g+a\bv_1+b\bv_2}{\sqrt{m^2-a^2-b^2}}$$
and is denoted by $F_{\mathcal{E}^{d}}^{d}(t,\theta)$. We call it {\it $\Delta_5$-dual surfaces} of $\mathcal{E}^{d}_{\g}$ under the assumption $\sigma<0$ and $m^2>a^2+b^2$.

 Assume $(a,b)\neq0$, under the Frenet type frame $\{\g,\bn_1,\bn_2,\m\}$, $F^{d}_{\mathcal{E}^{d}}(t,\theta)=\cosh\theta\m(t)+\sinh\theta(\frac{-M(t)\g(t)+A(t)\bn_1(t)}{\sqrt{M^2(t)-A^2(t)}}).$ Moreover, by equation (\ref{eq2.3}), we have
\begin{equation}\label{eqDE}
\begin{aligned}
(F^{d}_{\mathcal{E}^{d}})_{t}=&(\cosh\theta M+\sinh\theta\frac{(MA'-AM')(-A)}{(M^2-A^2)^{3/2}})\g+(-\cosh\theta A+\sinh\theta\frac{(MA'-AM')M}{(M^2-A^2)^{3/2}})\bn_1\\
&+(\frac{\sinh\theta AN}{\sqrt{M^2-A^2}})\bn_2
-\sinh\theta\sqrt{M^2-A^2}\m,\\
(F^{d}_{\mathcal{E}^{d}})_{\theta}=&\cosh\theta(\frac{-M\g+A\bn_1}{\sqrt{M^2-A^2}})-\sinh\theta\m.
\end{aligned}
\end{equation}
By the equations (\ref{eq4.26}) and (\ref{eqDE}), we can set  $\lambda_{\mathcal{E}^{d}}:I\times\Bbb R\rightarrow \Bbb R$ as follows:
\begin{equation}\label{eq5.38}
\lambda_{\mathcal{E}^{d}}(t,\theta)=\rm{ det}(F^{d}_{\mathcal{E}^{d}},(F^{d}_{\mathcal{E}^{d}})_{t},(F^{d}_{\mathcal{E}^{d}})_{\theta},\mathcal{E}^{d}_{\g})(t,\theta)=-\sinh\theta\frac{\sqrt{-\sigma_{F} (t)}}{M^2-A^2}(t).
\end{equation}

Then we have singular set $S(F^{d}_{\mathcal{E}^{d}})=\{(t,\theta)|\sinh\theta=0\}$. Moreover, by equation (\ref{eqDE}), we have $(F^{d}_{\mathcal{E}^{d}},\mathcal{E}^{d}_{\g})\in\Delta_5$ is isotropic. Thus $F^{d}_{\mathcal{E}^{d}}$ and $\mathcal{E}^{d}_{\g}$ are $\Delta_5$-dual each other.

\begin{Th}\label{Th5.4}
Under the above notation, we assume point $p_0=(t_0,\theta_0) \in S(F^{d}_{\mathcal{E}^{d}})$ under the assumption $\sigma_{F}<0$ and $M^2>A^2$. Then we have the followings.\vspace{2mm}\\
\mbox{}\quad
{\rm (1)} $p_0$ is a cuspidal edge if and only if
$$\frac{AN\sqrt{M^2-A^2}(MA'-AM')'-(AN\sqrt{M^2-A^2})'(MA'-AM')}{(MA'-M'A)^2+A^2N^2(M^2-A^2)}(t_0)
-\frac{MN}{\sqrt{M^2-A^2}}(t_0)\neq0,$$
\\
\mbox{}\quad
{\rm (2)} $p_0$ is cuspidal cross cap if and only if
$$\frac{AN\sqrt{M^2-A^2}(MA'-AM')'-(AN\sqrt{M^2-A^2})'(MA'-AM')}{(MA'-M'A)^2+A^2N^2(M^2-A^2)}(t_0)
-\frac{MN}{\sqrt{M^2-A^2}}(t_0)=0,$$
$$\frac{d}{dt}(\frac{AN\sqrt{M^2-A^2}(MA'-AM')'-(AN\sqrt{M^2-A^2})'(MA'-AM')}{(MA'-M'A)^2+A^2N^2(M^2-A^2)}
-\frac{MN}{\sqrt{M^2-A^2}})(t_0)\neq0.$$
\end{Th}
\proof
Since $(\lambda_{\mathcal{E}^{d}})_{\theta}\neq0$, by implicit theorem, then we obtain non-degenerate singular set $\xi_{\mathcal{E}^{d}}(t)=(t,\theta(t))$. And we set the null vector field $\eta_{\mathcal{E}^{d}}=\frac{\partial}{\partial t}+\sqrt{M^2-A^2}\frac{\partial}{\partial \theta}$ on $\xi_{\mathcal{E}^{d}}(t)$. Since $(\lambda_{\mathcal{E}^{d}})_{t}=0$, we obtain that  $(\xi_{\mathcal{E}^{d}})'(t)=(1,0)$. Then we consider
\begin{equation}
 {\rm det}(\eta_{\mathcal{E}^{d}},(\xi_{\mathcal{E}^{d}})^{'})(t)={\rm det}\begin{pmatrix}
1&\sqrt{M^2-A^2}\\
1&0\\
\end{pmatrix}(t)=-\sqrt{M^2-A^2}(t)\neq0.
\end{equation}
Thus, $p_0$ is a non-degenerate singular point of first kind. Furthermore, by a long but straightforward computation, we deduce that $\eta_{\mathcal{E}^{d}}\mathcal{E}_{\g}^{d}\neq0$ if and only if
\begin{equation}
\begin{aligned}
&\frac{AN\sqrt{M^2-A^2}(MA'-AM')'-(AN\sqrt{M^2-A^2})'(MA'-AM')}{(MA'-M'A)^2+A^2N^2(M^2-A^2)}
-\frac{MN}{\sqrt{M^2-A^2}}\neq0.\\
\end{aligned}
\end{equation}

The $\eta_{\mathcal{E}^{d}}\mathcal{E}_{\g}^{}\neq0$ proves that $F_{\mathcal{E}^{d}}^{d}$ is a front. Thus, the assertion (1) holds.

 Let us define a lift $\omega^{d}:U\rightarrow T^{*}_{F^{d}_{\mathcal{E}^{d}}}S^3_1$ by
$$\omega_{p}^{d}(\bv)=\langle \bv, \mathcal{E}_{\g}^{d}(p)\rangle, \bv\in  T_{F^{d}_{\mathcal{E}^{d}}}S^3_1, $$
where $p  $ is in the neighborhood $U$ of $F^{d}_{\mathcal{E}^{d}}$. We also define $\pi\circ\omega^{d}:U\rightarrow F^{d}_{\mathcal{E}^{d}}$, then $(\pi\circ\omega^{d})_{*}(r)=dF^{d}_{\mathcal{E}^{d}}(r)$ for any vector $r\in T_{p}U$. Since $\langle dF^{d}_{\mathcal{E}^{d}}(r),\mathcal{E}_{\g}^{d}(p)\rangle$=0, then we have $(\pi\circ\omega^{d})_{*}(T_{p}U)\subset\ker\omega_{p}^{d}$. This indicates that $\omega^{d}$ is the admissible lift of $F^{d}_{\mathcal{E}^{d}}$.
 Moreover, by a direct calculation, we obtain that $\frac{(\mathcal{E}_{\g}^{d})^{'}}{\|(\mathcal{E}_{\g}^{d})^{'}\|}$ is transversal to subspace $dF^{d}_{\mathcal{E}^{d}}(T_{p}U)$ at singular point $p$ and satisfies that $\langle \frac{(\mathcal{E}_{\g}^{d})^{'}}{\|(\mathcal{E}_{\g}^{d})^{'}\|}, \mathcal{E}_{\g}^{d}\rangle=0$. Putting \begin{equation*}
\begin{aligned}
\bX(t)=
\begin{cases}
\frac{(\mathcal{E}_{\g}^{d})^{'}(t)}{\|(\mathcal{E}_{\g}^{d})^{'}(t)\|}, \ \ &{\rm if} ~(\mathcal{E}_{\g}^{d})^{'}(t)\neq0,\\
{\lim\limits_{t\rightarrow +t_0}}\frac{(\mathcal{E}_{\g}^{d})^{'}(t)}{\|(\mathcal{E}_{\g}^{d})^{'}\|(t)}, \ \ &{\rm if}~ (\mathcal{E}_{\g}^{d})^{'}(t_0)=0,
\end{cases}
\end{aligned}
\end{equation*}
then $\bX$ ia a non-zero vector field along $F^{d}_{\mathcal{E}^{d}}|_{S(F^{d}_{\mathcal{E}^{d}})}$.
We set $\psi^{d}=\langle \eta_{\mathcal{E}^{d}}\bX, \mathcal{E}_{\g}^{d}\rangle$ on singular set $S(F^{d}_{\mathcal{E}^{d}})$.
 Since $\langle\bX, \mathcal{E}_{\g}^{d}\rangle=0$, then we have $\psi^{d}=-\langle \bX, \eta_{\mathcal{E}^{d}}\mathcal{E}^{d}_{\g}\rangle=-\langle \bX,(\mathcal{E}_{\g}^{d})^{'}\rangle$.
To make the following statement clear, we aseume that $(\mathcal{E}_{\g}^{d})^{'}(0)=0$. By the definition of $\bX$,
\begin{equation}
\begin{aligned}
\psi^{d}(t)=
\begin{cases}
-\langle\frac{(\mathcal{E}_{\g}^{d})^{'}(t)}{\|\mathcal{E}_{\g}^{d})^{'}\|(t)},(\mathcal{E}_{\g}^{d})^{'}\rangle, \ \ &{\rm if} ~(\mathcal{E}_{\g}^{d})^{'}(t)\neq0,\\
-\langle{\lim\limits_{t\rightarrow +0}}\frac{(\mathcal{E}_{\g}^{d})^{'}(t)}{\|\mathcal{E}_{\g}^{d})^{'}\|(t)},(\mathcal{E}_{\g}^{d})^{'}(0)\rangle \ \ &{\rm if}~ (\mathcal{E}_{\g}^{d})^{'}(0)=0
\end{cases}
\end{aligned}
\end{equation}
Since $(\mathcal{E}_{\g}^{d})^{'}(0)=0$, then we can find a deleted neighborhood of $t=0$ such that $(\mathcal{E}_{\g}^{d})^{'}(t)\neq0$. Moreover, we can deduce that $\psi^{d}(t)=-\frac{(\varepsilon^{d})^2(t)}{|\varepsilon^{d}(t)|}\neq0$ under $t\neq0$, where

\begin{equation*}
\begin{aligned}
&\varepsilon^{d}=\frac{AN\sqrt{M^2-A^2}(MA'-AM')'-(AN\sqrt{M^2-A^2})'(MA'-AM')}{(MA'-M'A)^2+A^2N^2(M^2-A^2)}\\
&-\frac{MN}{\sqrt{M^2-A^2}}.\\
\end{aligned}
\end{equation*}
satisfies $\varepsilon^{d}(t)=0$ if and only if $(\mathcal{E}_{\g}^{d})^{'}(t)=0$.
 Furthermore, we can assume that $\varepsilon^{d}(t)>0$ on the open interval $(0,t)$.
  By L'Hopital's rule, $\psi^{d}(0)=-{\lim\limits_{t\rightarrow +0}}\frac{(\varepsilon^{d})^2(t)}{|\varepsilon^{d}(t)|}=-2\varepsilon^{d}(0)=0$.
 Since $\frac{d}{dt}\psi^{d}(t)=-2(\varepsilon^{d})'(t)$ on the open interval $(0,t)$, then $\frac{d}{dt}\psi^{d}(0)=-{\lim\limits_{t\rightarrow +0}}2(\varepsilon^{d})'(t)=-2(\varepsilon^{d})'(0)$.

By the above, $ \psi^{d}(t)=0$ if and only if $t=0$ if and only if $\varepsilon^d(t)=0$. Thus, for a singular point $p_0=(t_0,\theta(t_0))$, $ \psi^{d}(t_0)=0$ and $ \frac{d}{dt}\psi^{d
}(t_0)\neq0$ if and only if $\varepsilon^{d}(t_0)=0$ and $(\varepsilon^{d})'(t_0)\neq0$. Then $ \psi^{d}(t_0)=0$ and $ \frac{d}{dt}\psi^{d}(t_0)\neq0$ prove that the assertion (2) holds.
\enD

Now we give the relationship between the $\Delta_5$-dual surfaces $F^{d}_{\mathcal{E}^{d}}$ and $\mathcal{E}^{d}_{\g}$.
\begin{pro}\label{pro5.5}
 We assume point $p_0=(t_0,\theta_0) \in S(F^{d}_{\mathcal{E}^{d}})$, then \vspace{2mm}\\
\mbox{}\quad
{\rm (1)} $p_0$ is a cuspidal edge if and only if  $\mathcal{E}^{d}_{\g}$ at $t_0$ is a regular point.\\
\mbox{}\quad
{\rm (2)} $p_0$ is cuspidal cross cap if and only if $\mathcal{E}^{d}_{\g}$ is {\rm (2,3,4)}-cusp at $t_0$
 \end{pro}
 \proof By the proof of Theorem \ref{Th5.4}, we can easily observe that $p_0$ is a cuspidal edge if and only if $(\mathcal{E}^{d}_{\g})'(t_0)\neq0$ and $p_0$ is a cuspidal cross cap if and only if $(\varepsilon^d)'(t_0)=0$ and $(\varepsilon^d)''(t_0)\neq0$. By a direct calculation, from the proof of Theorem \ref{Th4.8}, we obtain that $(\varepsilon^d)'(t_0)=0$ and $(\varepsilon^d)''(t_0)\neq0$ if and if $\mathcal{E}^{d}_{\g}$ is {\rm (2,3,4)}-cusp at $t_0$.
 \enD
By comparing the Theorem \ref{Th4.8} with Proposition \ref{pro5.5}, we have the following corollary.
 \begin{co}
  We assume point $p_0=(t_0,\theta_0)\in S(\mathcal{F}^{d})$, then there exists point $q_0=(t_0,\bar{\theta_0}) \in S(F^{d}_{\mathcal{E}^{d}})$ so that \vspace{2mm}\\
\mbox{}\quad
{\rm (1)} $\mathcal{F}^{d}$ at $(t_0,\theta_0)$ is cuspidal edge if and only if $F^{d}_{\mathcal{E}^{d}}$ at $(t_0,\bar{\theta_0})$ is a cuspidal edge\\
 \mbox{}\quad
  {\rm (2)} $\mathcal{F}^{d}$ at $(t_0,\theta_0)$ is swallowtail if and only if $F^{d}_{\mathcal{E}^{d}}$ at $(t_0,\bar{\theta_0})$ is a cuspidal cross cap
 \end{co}
\section{Dualities of conditions for the singular points}

In this section, we consider the self-dualities of cuspidal edge under the assumption $MA'-M'A\equiv0$.
\subsection{Dualities of singularities between $\mathcal{F}^{h}$ and $F_{\mathcal{E}^{h}}^{d}$}
In this subsection, we consider the dualities of singularities between $\mathcal{F}^{h}$ and $F_{\mathcal{E}^{h}}^{d}$ under the assumption $MA'-M'A\equiv0$. For the framed curve $(\g,\bn_1,\bn_2):I\rightarrow H^3\times \Delta_{5}$ with $A^2>M^2$, we obtain that $\mathcal{F}^{h}$ and $\m$ are $\Delta_1$-dual each other. By Theorem \ref{Th3.7}, we have that $\mathcal{E}_{\g}^{h}$ is the locus of non-degenerate singular values of $\mathcal{F}^{h}$ under $\sigma_{F}>0$. On the other hand, $F_{\mathcal{E}^{h}}^{d}$ and $\mathcal{E}_{\g}^{h}$ are $\Delta_1$-dual each other under $\sigma_{F}>0$.
By Theorem \ref{Th4.1}, we obtain that $\m$ is the locus of non-degenerate singular values of $F_{\mathcal{E}^{h}}^{d}$. Since we deduce that $A^2>M^2$ by $\sigma_{F}>0$, then we have the following diagram.
$$
\xymatrixcolsep{10pc}\xymatrix{
\mathcal{F}^{h} \ar@{<->}[d]_{\Delta_1-dual}^{A^2>M^2}  \ar@{->}[r]^{taking~singular~value}_{\sigma_{F}>0} &\mathcal{E}_{\g}^{h}\ar@{<->}[d]^{\Delta_1-dual}_{\sigma_{F}>0} \\
\m& \ar[l]^{taking~singular~value}    F_{\mathcal{E}^{h}}^{d}}
$$
In order to investigate the dual relations of conditions for singularities on dual surfaces, we assume that $MA'-M'A\equiv0$ so that $\{(t,0)\}$ on $\mathcal{F}^{h}$ is always singular.
 Furthermore, we give the following classifications of singular points of $\mathcal{F}^{h}$ under the assumption $MA'-M'A\equiv0$.
\begin{Th}\label{Th4.2}
Under the assumption $MA'-M'A\equiv0$, $p_0\in S(\mathcal{F}^{h})$. Then we have the followings.\vspace{2mm}\\
\noindent
{\rm (a)} If $N(t_0)\neq0$, then\\
\mbox{}\quad
{\rm (1)} $p_0$ is cuspidal edge if and if $M(t_0)\neq0$.\\
\mbox{}\quad
{\rm (2)} $p_0$ is never swallowtail and cuspidal cross cap. \\
\noindent
{\rm (b)}  If $N(t_0)=0$, then\\
\mbox{}\quad
{\rm (1)} $p_0$ is cuspidal edge if and if $\theta_0\neq0$, $(AN\sqrt{A^2-M^2})'(t_0)\neq0$.\\
\mbox{}\quad
{\rm (2)} $p_0$ is never cuspidal beaks, cuspidal lips, swallowtail and cuspidal cross cap.
\end{Th}
Under the $MA'-M'A\equiv0$, by (\ref{eq4.18}), we can see that $(t,0)$ on $F_{\mathcal{E}^{h}}^{d}$ is always singular. Furthermore, we calculate that $\sigma_{F}>0$ if and only if $A^2>M^2$ and $N\neq0$ under the $MA'-M'A\equiv0$. Then we also give the following classifications of singular points of $F_{\mathcal{E}^{h}}^{d}$ under the assumption $MA'-M'A\equiv0$, $A^2>M^2$ and $N\neq0$.
\begin{Th}\label{Th4.3}
Under the assumption $MA'-M'A\equiv0$, $A^2>M^2$ and $N\neq0$, $p_0\in S(F_{\mathcal{E}^{h}}^{d})$. Then we have the followings.\vspace{2mm}\\
\mbox{}\quad
{\rm (a)} $p_0$ is cuspidal edge if and if $M(t_0)\neq0$.\\
\mbox{}\quad
{\rm (b)} $p_0$ is never cuspidal cross cap.
\end{Th}
\begin{re}
By the proofs of Theorem \ref{Th3.3} and Theorem \ref{Th4.1}, we can simply obtain the proofs of Theorem \ref{Th4.2} and Theorem \ref{Th4.3}. Thus, we omit them.
 \end{re}
 By comparing Theorem \ref{Th4.2} and Theorem \ref{Th4.3}, we summarize the conditions on the following Table.
 \begin{equation*}
\begin{tabular}{|c||c||c||c||c|}
\hline
surfaces&duality&S=\{(t,0)\}&CE\\
\hline
$\mathcal{F}^{h}$&$A^{2}>M^{2}$&alwarys&$ N(t_0)\neq0$$M(t_0)\neq0$\\
\hline
$F_{\mathcal{E}^{h}}^{d}$&$A^{2}>M^{2}$, $N\neq0$&alwarys&$M(t_0)\neq0$\\
\hline
\end{tabular}
\end{equation*}
From the above table, we can easily observe the self-duality of cuspidal edge under $MA'-M'A\equiv0$.

\subsection{Dualities of singularities between $\mathcal{F}^{d}$ and $F_{\mathcal{E}^{d}}^{d}$}
In this subsection, we consider the dualities of singularities between $\mathcal{F}^{d}$ and $F_{\mathcal{E}^{d}}^{d}$ under the assumption $MA'-M'A\equiv0$.
For the framed curve $(\g,\bn_1,\bn_2):I\rightarrow H^3\times \Delta_{5}$ with $M^2>A^2$, we obtain that $\mathcal{F}^{d}$ and $\m$ are $\Delta_5$-dual each other. By Theorem \ref{Th4.8}, we deduce that $\mathcal{E}_{\g}^{d}$ is the locus of non-degenerate singular values of $\mathcal{F}^{d}$ under $\sigma_{F}<0$. On the other hand, $F_{\mathcal{E}^{d}}^{d}$ and $\mathcal{E}_{\g}^{d}$ are $\Delta_5$-dual each other under $\sigma_{F}<0$ and $M^2>A^2$. By Theorem \ref{Th5.4}, we obtain that $\m$ is the locus of non-degenerate singular values of $F_{\mathcal{E}^{d}}^{d}$. Then we have the following diagram.

$$
\xymatrixcolsep{10pc}\xymatrix{
\mathcal{F}^{d} \ar@{<->}[d]_{\Delta_5-dual}^{M^2>A^2}  \ar@{->}[r]^{taking~singular~value}_{\sigma_{F}<0} &\mathcal{E}_{\g}^{d}\ar@{<->}[d]^{\Delta_5-dual}_{\sigma_{F}<0,M^2>A^2} \\
\m& \ar[l]^{taking~singular~value}    F_{\mathcal{E}^{d}}^{d}}
$$

We also assume that $MA'-M'A\equiv0$ so that $\{(t,0)\}$ on $\mathcal{F}^{d}$ is always singular.
 Furthermore, we give the following classifications of singular points of $F^{d}$ under the assumption $MA'-M'A\equiv0$.
\begin{Th}\label{Th6.4}
Under the assumption $MA'-M'A\equiv0$, $p_0\in S(\mathcal{F}^{d})$. Then we have the followings.\vspace{2mm}\\
\noindent
{\rm (a)} If $N(t_0)\neq0$, then\\
\mbox{}\quad
{\rm (1)} $p_0$ is cuspidal edge.\\
\mbox{}\quad
{\rm (2)} $p_0$ is never swallowtail and cuspidal cross cap. \\
\noindent
{\rm (b)}  If $N(t_0)=0$, then\\
\mbox{}\quad
{\rm (1)} $p_0$ is cuspidal edge if and if $\theta_0\neq0$ or $\pi$, $(AN\sqrt{M^2-A^2})'(t_0)\neq0$.\\
\mbox{}\quad
{\rm (2)} $p_0$ is never cuspidal beaks, cuspidal lips, swallowtail and cuspidal cross cap.
\end{Th}

Under the $MA'-M'A\equiv0$, by (\ref{eq5.38}), we can see that $(t,0)$ on $F_{\mathcal{E}^{d}}^{d}$ is always singular. Furthermore, we calculate that $\sigma_{F}<0$, $M^2>A^2$ if and only if $M^2>A^2$ and $N\neq0$ under the $MA'-M'A\equiv0$. Then we also give the following classifications of singular points of $F_{\mathcal{E}^{d}}^{d}$ under the assumption $MA'-M'A\equiv0$, $M^2>A^2$ and $N\neq0$.

\begin{Th}\label{Th6.5}
Under the assumption $MA'-M'A\equiv0$, $M^2>A^2$ and $N\neq0$, $p_0\in S(F_{\mathcal{E}^{d}}^{d})$. Then we have the followings.\vspace{2mm}\\
\mbox{}\quad
{\rm (a)} $p_0$ is cuspidal edge.\\
\mbox{}\quad
{\rm (b)} $p_0$ is never cuspidal cross cap.
\end{Th}
\begin{re}
By the proofs of Theorem \ref{Th4.6} and Theorem \ref{Th5.4}, we simply obtain the proofs of Theorem \ref{Th6.4} and Theorem \ref{Th6.5}. Thus, we omit them.
 \end{re}
By comparing Theorem \ref{Th6.4} and Theorem \ref{Th6.5}, we summarize the conditions on the following Table.
 \begin{equation*}
\begin{tabular}{|c||c||c||c||c|}
\hline
surfaces&duality&S=\{(t,0)\}&CE\\
\hline
$\mathcal{F}^{d}$&$M^{2}>A^{2}$&alwarys&$N(t_0)\neq0$\\
\hline
$F_{\mathcal{E}^{d}}^{d}$&$M^{2}>A^{2}$, $N\neq0$&alwarys&alwarys\\
\hline
\end{tabular}
\end{equation*}
From the above table, we can easily observe the self-duality of cuspidal edge under $MA'-M'A\equiv0$.

\end{sloppypar}
\typeout{get arXiv to do 4 passes: Label(s) may have changed. Rerun}

\begin{thebibliography}{99999}
{\renewcommand{\baselinestretch}{1} \small


\bibitem{1}
Arnol'd V. I., Gusein-Zade S. M., Varchenko A. N., Singularities of differentiable maps, vol. I. Boston: Birkh\"{a}user, 1985

\bibitem{2}
Bruce, J. W., Giblin, P. J., Curves and singularities, 2nd; Cambridge University press: Cambridge, UK, 1992.

\bibitem{3}
Chen L., Izumiya S. A mandala of Legendrian dualities for pseudo-spheres in semi-Euclidean space. Proc Japan Acad Ser A Math Sci, 2009, 85 : 49-54

\bibitem{4}
Fuchs, D. Evolutes and involutes of spatial curves. Amer. Math. Monthly (2013) 120, 217-231.

\bibitem{5}
Fujimori S, Saji K, Umehara M, et al. Singularities of maximal surfaces. Math Z, 2008, 259: 827-848

\bibitem{6}
Fukunaga T., Tahakashi M., Evolutes of front in the Euclidean Plane, J. Singul. 10 (2014) 92-107.

\bibitem{7}
Fukunaga T., Tahakashi M., Evolutes and Involutes of Frontals in the Euclidean Plane, Demonstratio Math. 48 (2015) 147-166.


\bibitem{8}
 Gray A., Abbena E., Salamon S., Modern differential geometry of curves and surfaces with Mathematica, 3rd; Chapman and Hall/CRC: Boca Raton, FL, USA, 2006.

\bibitem{9}
Honda S., Takahashi M., Framed curves in the Euclidean space. Adv. Geom. 16 (2016) 265-276

\bibitem{10}
Honda, S., Takahashi, M., Evolutes and focal surfaces of framed immersions in the Euclidean space. Proc. Roy. Soc. Edinburgh Sect. A (2020) 150, 497-516.

\bibitem{11}
 Porteous, I. R. Geometric differentiaton for the intelligence of curves and surfaces, 2nd; Cambridge University Press: Cambridge, UK, 2001.

\bibitem{12}
 Porteous I. R., Some remarks on duality in S3, Geometry and Topology of Caustics, Banach Center
Publications 50, Polish Academy of Sciences, Warsaw, 2004, 217¨C226.

\bibitem{13}
Hayashi R., Izumiya S., Sato T., Focal surfaces and evolutes of curves in hyperbolic space. Commun. Korean Math. Soc. 32 (2017), no.1, 147-163.


\bibitem{14}
Izumiya S., Saji K., The mandala of Legendrian dualities for pseudo-spheres in Lorentz-Minkowski space and ``flat" spacelike surfaces. J Singul, 2010, 2: 92-127


\bibitem{15}
Izumiya S., Legendrian dualities and spacelike hypersurfaces in the lightcone. Mosc Math J, 2009, 9: 325-357

\bibitem{16}
Izumiya S., Jiang Y., and Pei D. H., Lightcone dualities for curves in the sphere, Q. J. Math. 64, 221¨C234 (2013)

\bibitem{17}
Izumiya S., Jiang Y., and Sato T., Lightcone dualities for curves in the 3-sphere, J. Math. Phys. 54, 063511 (2013)


\bibitem{18}
Kokubu M., Rossman W., Saji K., et al. Singularities of flat fronts in hyperbolic 3-space. Pacific J Math, 2005, 221: 303-351

\bibitem{19}
Martins L. F., Saji K., Umehara M., et al. Behavior of Gaussian curvature and mean curvature near non-degenerate singular points on wave fronts. Springer Proc Math Stat, 2016, 154: 247-281.




\bibitem{20}
Saji K., Umehara M., Yamada K., The geometry of fronts. Ann of Math (2), 2009, 169: 491-529


\bibitem{21}
Saji K., Umehara M., Yamada K., $A_{k}$ singularities of wave fronts. Math Proc Cambridge Philos Soc, 2009, 146: 731-746




\bibitem{22}
Zhou A., Yao K., Pei D. H., {\it k}-type hyperbolic framed slant helices in hyperbolic 3-space. Filomat 38 (2024), no.1, 3839-3850.






}
\end{thebibliography}
\end{document}